\documentclass[10pt]{article}
\usepackage{graphicx,epsfig}
\usepackage{amsmath}

 \usepackage{color}      % use if color is used in text

\usepackage{url}

\usepackage{tabu}

\usepackage{array}
\usepackage{makecell}

\usepackage{multirow}

\usepackage[thinlines]{easytable}

\usepackage{placeins}
\usepackage{indentfirst}
\usepackage{multirow}
\usepackage{url}
\usepackage{enumerate}
\usepackage{graphicx}

\usepackage{booktabs}
 \newcommand{\ra}[1]{\renewcommand{\arraystretch}{#1}}
 
\usepackage{subcaption}
\usepackage{mathtools}

\usepackage{subcaption}
\usepackage[driverfallback=dvipdfm,
bookmarks=true,    %show bookmarks bar?
bookmarksopen=true,
bookmarksopenlevel=\maxdimen,
bookmarksdepth=10
]{hyperref}
\hypersetup{
pdfcenterwindow=true, 
pdfpagelabels=true,
pagebackref=true,            
hypertexnames=true,
plainpages=false,
unicode=false,     %non-Latin characters in Acrobat$B!G(Bs bookmarks
pdftoolbar=true,                % show Acrobat$B!G(Bs toolbar?
    pdfmenubar=true,          % show Acrobat$B!G(Bs menu?
    pdffitwindow=false,         % window fit to page when opened
    pdfstartview={FitH},        % fits the width of the page to the window
    pdftitle={An Efficient Implementation of Edge-Based Discretization without Forming Dual Control Volumes},
    pdfauthor={Hiroaki Nishikawa},    % author
    pdfsubject={An Efficient Implementation of Edge-Based Discretization without Forming Dual Control Volumes},       % subject of the document
    pdfcreator={Hiroaki Nishikawa},   % creator of the document
    pdfproducer={Hiroaki Nishikawa}, % producer of the document
    baseurl={https://www.researchgate.net/profile/Hiroaki-Nishikawa-2},
    pdfkeywords={}, % list of keywords
    pdfnewwindow=true,     % links in new window
    colorlinks=true,       % false: boxed links; true: colored links
    linkcolor=black,       % color of internal links
    citecolor=black,       % color of links to bibliography
    filecolor=black,       % color of file links
    urlcolor=black,        % color of external links
  breaklinks=true,
  hyperfigures=true,
  backref=true,
breaklinks=false, 
pdfpagelabels,      %%%%% Options here and below are needed to be compatible with
pagebackref,        %%%%% \numberwithin command.
hypertexnames=true, %%%%%
plainpages=false,   %%%%%
naturalnames        %%%%%
}
\usepackage[all]{hypcap}
\numberwithin{equation}{section}
\numberwithin{subsubsection}{subsection}
\numberwithin{subsection}{section}

\usepackage{siunitx}

\usepackage{cancel}

\title{\bf An Efficient Implementation of Edge-Based Discretization without Forming Dual Control Volumes} %
\author{ 
{Hiroaki Nishikawa}\thanks{Research Fellow ({hiro@nianet.org})  }\\
  {\normalsize\itshape National Institute of Aerospace, Hampton, VA 23666, USA}
}

\def\o6{\frac{1}{6}} 
\usepackage{cleveref}

\begin{document}

\maketitle

\begin{abstract}
This paper shows that lumped directed-area vectors at edges and dual control volumes required to implement the edge-based discretization can be computed without explicitly defining the dual control volume around each node for triangular and tetrahedral grids. It is a simpler implementation because there is no need to form a dual control volume by connecting edge-midpoints, face centroids, and element centroids, and also reduces the time for computing lumped directed-area vectors for a given grid, especially for tetrahedral grids. The speed-up achieved by the proposed algorithm may not be large enough to greatly impact the overall simulation time, but the proposed algorithm is expected to serve as a major stepping stone towards extending the edge-based discretization to four dimensions and beyond (e.g., space-time simulations). Efficient algorithms for computing lumped directed-area vectors and dual volumes without forming dual volumes are presented, and their implementations are described and compared with traditional algorithms in terms of complexity as well as actual computing time for a given grid. 
\end{abstract}

%%%%%%%%%%%%%%%%%%%%%%%%%%%%%%%%%%%%%%%%%%%%%%%%
\section{Introduction}
\label{intro}
%%%%%%%%%%%%%%%%%%%%%%%%%%%%%%%%%%%%%%%%%%%%%%%%

\indent
 
The node-centered edge-based discretization method is widely used today in practical unstructured-grid computational fluid dynamics (CFD) solvers 
\cite{barth_AIAA1991,fun3d_manual:NASATM2016,nakashima_watanabe_nishikawa:Japan2014,dlr-tau-digital-x,mavriplis_long:AIAA2010,Luo_Baum_Lohner:AIAA2004-1103,KozubskayaAbalakinDervieux:AIAA2009,Haselbacher_Blazek_AIAAJ2000,sierra-primo:AIAA2002,Eliasson_EDGE:2001,fezoui_stoufflet:JCP1989,GaoHabashiFossatiIsolaBaruzzi_AIAA2017-0085}. It is a second-order accurate conservative discretization on arbitrary simplex-element (triangular or tetrahedral) grids with point-valued solutions stored at nodes \cite{Nishikawa_FSR:2020} that can be implemented in an efficient edge loop, where only solutions and gradients at two end nodes are needed to compute a numerical flux at each edge. The method has also been shown to achieve third-order accuracy with a single numerical flux per edge on arbitrary triangular \cite{Katz_Sankaran_JCP:2011,katz_sankaran:JSC_DOI,diskin_thomas:AIAA2012-0609} and 
tetrahedral grids \cite{liu_nishikawa_aiaa2016-3969,liu_nishikawa_aiaa2017-0081} without high-order grids \cite{liu_nishikawa_aiaa2016-3969,nishikawa_boundary_formula:JCP2015} nor second derivatives \cite{nishikawa_liu_source_quadrature:jcp2017}. Because of these special properties, it has been a subject for active CFD algorithm research towards automated CFD simulations with anisotropic viscous grid adaptation that can be performed efficiently with simplex-element grids \cite{Kleb_etal_aiaa2019-2948,ThompsonNishikawaPadway_aiaa_scitech2023}. 

The edge-based discretization method constructs a residual over a medial dual control volume around each node, which is formed by connecting edge midpoints, geometric centroids of elements, and geometric centroids of faces (in three dimensions). To implement the edge-based discretization method, we need to compute volumes of the medial dual control volumes and lumped directed-area vectors at all edges for a give grid. The lumped directed-area vector is defined at each edge as the sum of the directed-area vectors of dual-volume faces attached to that edge. Typically, these two grid metrics are computed in a loop over elements, accumulating local contributions at nodes for dual volumes and at edges for lumped directed-area vectors. This is done by partially forming a median dual control volume within each element. Therefore, the construction of the median dual control volume is essential. In extending the edge-based discretization method to four and higher dimensions (e.g., space-time simulations), however, we encounter a difficulty that it is not straightforward, although perhaps possible, to construct a median dual control volume and compute such local contributions in a higher-dimensional element such as a pentatope in four dimensions. Although various discretization methods have been extended to four-dimensional grids \cite{behr:IJNMF2008,MDG_aiaa2019-0642,Caplan:CAD2020,Foundation_STFEM:ANM2021,McCaughtryWatsonGeronKim:Aviation2022}, to the author's knowledge at the time of writing this paper, the edge-based discretization has not been applied to four- or higher-dimensional simplex-element grids. It motivated the author to develop an alternative implementation that does not require the construction of median dual control volumes, which is the main subject of this paper. Such an implementation is expected not only to greatly simplify the coding of an edge-based solver for higher dimensions, but also simplify and potentially improve the efficiency of edge-based solvers in two and three dimensions. Leaving the former as important future work, we focus here on the basic idea, and derive and demonstrate efficient computations of lumped directed-area vectors and dual volumes without forming median dual control volumes for triangular and tetrahedral grids.

The main idea is to express directed-area vectors of dual faces entirely in terms of directed-area vectors of element faces. In doing so, we would have to introduce extra vectors on side faces, but they all cancel out after summing local contributions over a set of elements sharing an edge of interest, for all interior edges. As a result, we obtain an efficient formula for computing lumped directed-area vectors at interior edges without explicitly constructing median dual control volumes. At boundary edges, the lumped directed-area vectors need to be corrected with extra vectors, but this correction can be performed efficiently in a loop over boundary elements, again, without explicitly forming dual volumes and dual faces. Moreover, dual volumes can also be computed efficiently, without forming median dual control volumes, in a loop over edges using the lumped directed-area vectors. Therefore, there is no need at all to form the median dual control volumes. It is possible that some of the proposed implementations have already known to some researchers especially in early years when the edge-based discretization was discussed in relation to a continuous Galerkin finite-element discretization in the 1990s (see, e.g., Ref.~\cite{barth_AIAA1991}). However, such implementations might not have been considered particularly attractive partly because the grid metric computation takes only a fraction of an entire CFD simulation time in many cases and the construction of dual volumes is not difficult in two and three dimensions, and also because extensions to four and higher dimensions have not been received much attention. In this paper, we focus on two- and three-dimensional simplex-element grids. Extensions and applications to higher dimensions will be reported in a separate paper. 

The paper is organized as follows. 
In Section 2, we describe the edge-based discretization. 
In Section 3, we begin with a traditional algorithm for computing lumped directed-area vectors based on median dual volumes and then introduce a more efficient algorithm that does not require median dual volumes. 
In Section 4, we discuss two approaches to computing median dual volumes without explicitly constructing median dual control volumes. 
In Section 5, we present numerical results comparing the efficiency of the new algorithms to traditional algorithms. 
In Section 6, we end the paper with concluding remarks.

% \cite{nishikawa_aiaa2017-4295}   Zero/nega vols
    
%%%%%%%%%%%%%%%%%%%%%%%%%%%%%%%%%%%%%%%%%%%%%%%%
\section{Edge-Based Discretization}
\label{eb_discretization}
%%%%%%%%%%%%%%%%%%%%%%%%%%%%%%%%%%%%%%%%%%%%%%%%

\indent

%where $T$ is the temperature. As mentioned earlier, the proposed methodology will be generally applicable to other hyperbolic systems including the hyperbolic diffusion system \cite{nishikawa_onetwothree_diffusion:JCP2014} and the hyperbolic Navier-Stokes systems %\cite{liu_nishikawa_aiaa2016-3969,NakashimaWatanabeNishikawa_AIAA2016-1101,LiLouNishikawaLuo_HNSrDG:JCP2021}. 

\begin{figure}[htbp!]
    \centering
          \begin{subfigure}[t]{0.48\textwidth}
        \includegraphics[width=0.9\textwidth,trim=0 0 0 0,clip]{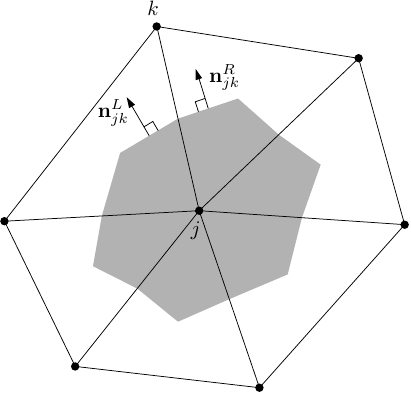}
          \caption{Triangular grid.}
          \label{fig:dual_vol_2d}
      \end{subfigure}
      \hfill
          \begin{subfigure}[t]{0.48\textwidth}
        \includegraphics[width=1.2\textwidth,trim=2 2 2 2,clip]{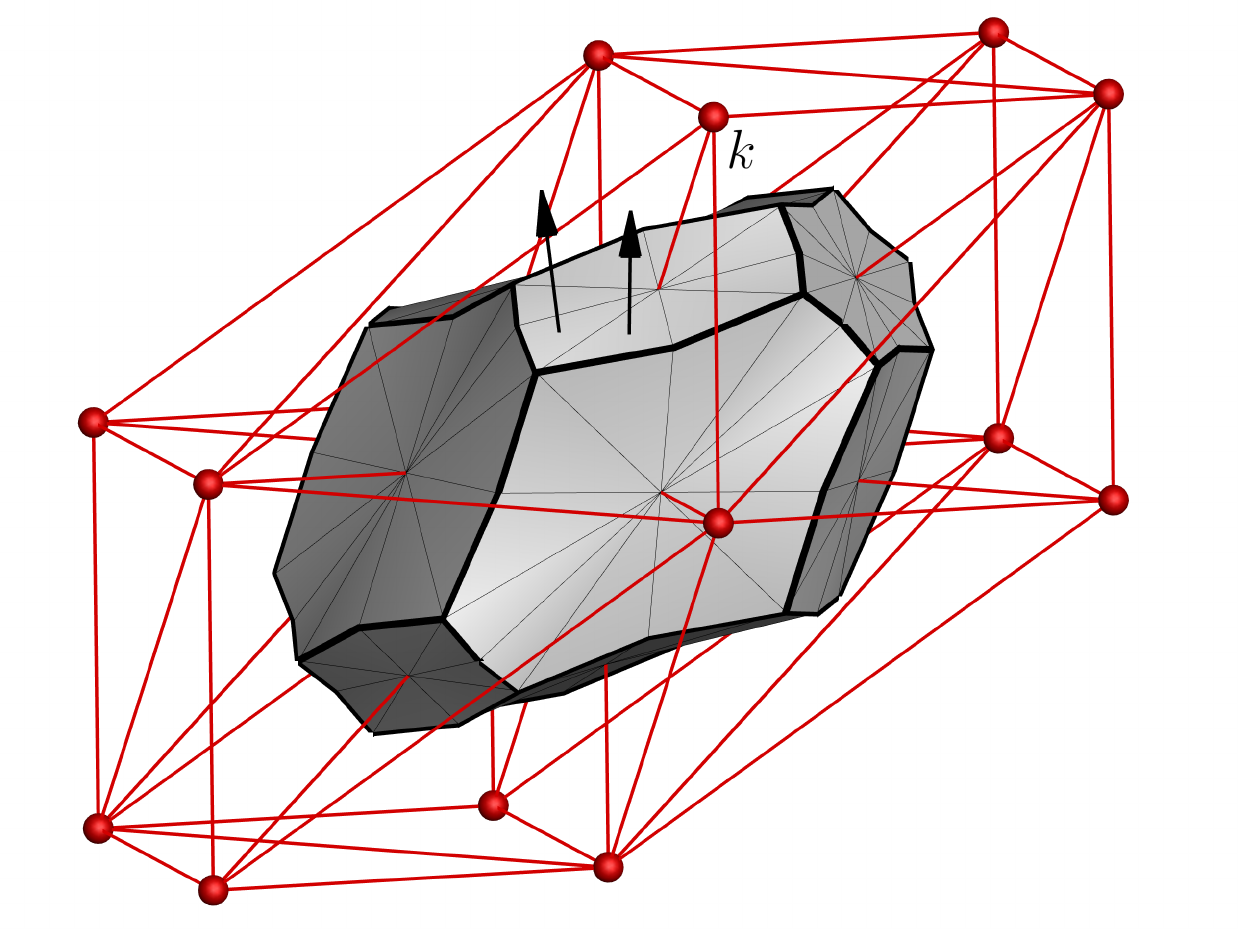}
          \caption{Tetrahedral grid.}
          \label{fig:dual_vol_3d}
      \end{subfigure}
            \caption{
\label{fig:njk_contribution}%
Median dual control volumes in two and three dimensions, constructed with edge midpoints and geometric centroids of triangular elements in two dimensions, and edge midpoints, geometric centroids of tetrahedral elements, and geometric centroids of triangular faces in three dimensions. Note that the center node $j$ in the tetrahedral grid is located inside the median control volume and is not seen in the figure. } 
\end{figure}

Consider a conservation law in the differential form:
\begin{eqnarray}
\frac{ \partial {\bf u} }{\partial t}  + \mbox{div}{\cal F} = {\bf s}, 
\label{diff_form}
\end{eqnarray}
where $t$ is time, ${\bf u}$ is a vector of solution variables, ${\cal F}$ is a flux tensor, and ${\bf s}$ is a source (or forcing) term vector. 
The edge-based discretization method is a method of discretizing a conservation law (\ref{diff_form}) at a node in a conservative manner over a median dual control volume defined by connecting edge midpoints and geometric centroids of elements (and geometric centroids of element faces in three dimensions), with point-valued solutions stored at nodes. See Figures \ref{fig:dual_vol_2d} and \ref{fig:dual_vol_3d}. In this paper, we focus on simplex-element grids (triangular grids in two dimensions and tetrahedral grid in three dimensions): the edge-based discretization at a node $j$ is given by 
\begin{eqnarray}
 V_j \frac{ d {\bf u}_j }{d t}  +  \sum_{k \in \{ k_j\} } {\Phi}_{jk} | {\bf n}_{jk} | = {\bf s}_j \, V_j,
\label{threed_fv_semidiscrete_system_00}
\end{eqnarray}
where $V_j$ is the measure of the median dual control volume around the node $j$, $\{ k_j\}$ is a set of edge-connected neighbor nodes of the node $j$, ${\Phi}_{jk}$ is a numerical flux, and ${\bf n}_{jk}$ is the lumped directed-area vector, which is the sum of the directed-area vectors corresponding to the dual faces associated with all elements sharing the edge $\{ j, k \}$. In two dimensions, the lumped directed-area vector is defined as the sum of two directed-area vectors ${\bf n}_{jk}^L$ and ${\bf n}_{jk}^R$ as illustrated in Figure \ref{fig:dual_vol_2d}. In three dimensions, it is defined as the sum of the surface directed-area vectors of all the triangular dual faces associated with the edge $\{ j, k \}$ as illustrated in Figure \ref{fig:dual_vol_3d} (arrows are illustrated only for two of the eight relevant dual triangular faces). Note that each edge is oriented from one node to the other and the lumped directed-area vector is directed consistently. A simple way of defining the orientation is to define an edge as oriented from a node with a smaller node number to the other: ${\bf n}_{jk} $ pointing from $j$ to $k$, where $ j < k$. Throughout this paper, we assume that all edges are oriented in this way. 

%%%%%%%%%%%%%%%%%%%%%%%%%%%%%%%%%%%%%%%%%
%  
%%%%%%%%%%%%%%%%%%%%%%%%%%%%%%%%%%%%%%%%%
%%%%%%%%%%%%%%%%%%%%%%%%%%%%%%%%%%%%%%%%%
% 
%%%%%%%%%%%%%%%%%%%%%%%%%%%%%%%%%%%%%%%%%

The numerical flux is computed at the edge midpoint, for example, by an upwind flux,
\begin{eqnarray}
 {\Phi}_{jk} ( {\bf u}_L, {\bf u}_R,  \hat{\bf n}_{jk} ) = 
  \frac{1}{2} \left[     {\bf f}({\bf u}_L) + {\bf f}({\bf u}_R)      \right]   
    - \frac{1}{2} \left|   \frac{\partial {\bf f} }{\partial {\bf u}} \right| \left(  {\bf u}_R - {\bf u}_L \right ) ,  
  \label{upwind_flux}
\end{eqnarray}
where %${\bf w} = (\rho , u, v, w, p)$, 
${\bf f} = {\cal F} \cdot  \hat{\bf n}_{jk}$, $ \hat{\bf n}_{jk} =  {\bf n}_{jk}  /  | {\bf n}_{jk} |  $, and the subscripts $L$ and $R$ indicate values at the left and right sides of the edge midpoint. The left and right solution values, ${\bf u}_L$ and ${\bf u}_R$, are reconstructed from the solutions stored at $j$ and $k$, respectively, by an algorithm that is exact at least for linear functions: e.g., the U-MUSCL scheme \cite{burg_umuscl:AIAA2005-4999,PadwayNishikawa:AIAAJ2022} and the V4 scheme \cite{DebiezDerieux:CF2000}. The edge-based discretization described above is known to produce second-order accurate point-valued solutions at nodes on arbitrary simplex-element grids, and it can also produce third-order accurate point-valued solutions at nodes with minor modifications, as mentioned earlier. 

The flux balance $ \sum_{k \in \{ k_j\} } {\Phi}_{jk} | {\bf n}_{jk} |$ in the residual (\ref{threed_fv_semidiscrete_system_00}) can be computed by looping over edges: $e = 1, 2 , 3, \cdots, N_e$, where $N_e$ is the number of edges in a grid. More precisely, we first initialize the residual array ${\bf Res}_j = 0$ for all nodes, $k=1,2,3, \cdots, N_v$, where $N_v$ is the number of nodes, and then looping over edges, we compute the flux and distribute it to the left and right nodes at each edge $e$: 
\begin{eqnarray} 
 {\bf Res}_j = {\bf Res}_j +  {\Phi}_{jk}, \quad 
 {\bf Res}_k = {\bf Res}_k -  {\Phi}_{jk},
\end{eqnarray}
where $j$ and $k$ are the nodes of the edge $e$, which may be stored in a two-dimensional array, $edge(e,1) = j $ and $edge(e,2) = k$, and ${\bf n}_{jk}$ is defined as pointing from $j$ and $k$. If boundary conditions are imposed weakly, the residuals at boundary nodes need to be closed with boundary flux contributions, and it can be performed by a loop over boundary elements (see Ref.~\cite{nishikawa:AIAA2010}). If necessary, the time derivative and the source term can be added afterwards in a loop over nodes. 

%It can achieve third-order accuracy with minor modifications. See Refs.~\cite{Katz_Sankaran_JCP:2011,katz_sankaran:JSC_DOI,diskin_thomas:AIAA2012-0609,nishikawa_boundary_quadrature:JCP2015,nishikawa_liu_source_quadrature:jcp2017} for details. To implement boundary conditions in a weak manner, we need to close the residual at a boundary node and special boundary flux quadrature formulas are necessary, for both second- and third-order schemes, to preserve the design order of accuracy. See Appendix B of Ref.~\cite{nishikawa:AIAA2010} for the second-order scheme and Ref.~\cite{nishikawa_boundary_quadrature:JCP2015}. 

The edge-based discretization is efficient because it requires only a single numerical flux per edge even for third-order accuracy. This exceptional efficiency is due to the special property of the edge-based discretization defined with the lumped directed-area vectors and the median dual volume: it gives $ \sum_{k \in \{ k_j\} } {\Phi}_{jk} | {\bf n}_{jk} | / V_j  = ( \mbox{div} {\cal F} )_j$ exactly at a node $j$ on arbitrary simple-element grids for a linearly-varying (in space) flux ${\cal F}$, and the same is true for a quadratically-varying flux in the case of the third-order edge-based discretization \cite{nishikawa_liu_source_quadrature:jcp2017}.

Note that only two grid metrics need to be precomputed for a given grid to implement the second- and third-order edge-based discretizations: the lumped directed-area vectors at edges and the volumes of median dual volumes at nodes. Our main interest here is in the computations of these quantities. In the next section, we begin with algorithms to compute the lumped directed-area vectors, and then discuss algorithms for computing the median dual volumes. 

%Subsequently, we will derive new efficient boundary closure formulas from the efficient algorithm for the directed-area vectors. 

%Notice that the discretization is completely defined by the numerical flux, the dual volume, and the directed-area vectors. 
%Our interest here is in the dual volume and the directed-area vectors, which can be pre-computed and stored for a given grid, 
%and in deriving and identifying efficient formulas for computing them. 
%Although it will only reduce the computational time for the pre-computing step, it can be 

%%%%%%%%%%%%%%%%%%%%%%%%%%%%%%%%%%%%%%%%%%%%%%%%
%%%%%%%%%%%%%%%%%%%%%%%%%%%%%%%%%%%%%%%%%%%%%%%%
%%%%%%%%%%%%%%%%%%%%%%%%%%%%%%%%%%%%%%%%%%%%%%%%
%%%%%%%%%%%%%%%%%%%%%%%%%%%%%%%%%%%%%%%%%%%%%%%%
%%%%%%%%%%%%%%%%%%%%%%%%%%%%%%%%%%%%%%%%%%%%%%%%
\section{Efficient Computations of Lumped Directed-Area Vectors}
\label{eb_njk_vj}
%%%%%%%%%%%%%%%%%%%%%%%%%%%%%%%%%%%%%%%%%%%%%%%%
%%%%%%%%%%%%%%%%%%%%%%%%%%%%%%%%%%%%%%%%%%%%%%%%
%%%%%%%%%%%%%%%%%%%%%%%%%%%%%%%%%%%%%%%%%%%%%%%%
%%%%%%%%%%%%%%%%%%%%%%%%%%%%%%%%%%%%%%%%%%%%%%%%
%%%%%%%%%%%%%%%%%%%%%%%%%%%%%%%%%%%%%%%%%%%%%%%%

 \indent 
 
In this section, we first review a traditional algorithm for computing the lumped directed-area vectors for triangular and tetrahedral grids, and describe a more efficient algorithm proposed in this paper. Then, we provide a complexity analysis and discuss the relative cost of the two algorithms. 
%Here, we focus on interior edges, and discuss boundary edges in a later section.

%%%%%%%%%%%%%%%%%%%%%%%%%%%%%%%%%%%% 
\subsection{A Traditional Algorithm for Computing Lumped Directed-Area Vectors}
\label{eb_njk_typical}
%%%%%%%%%%%%%%%%%%%%%%%%%%%%%%%%%%%% 

%%%%%%%%%%%%%%%%%%%%%%%%%%%%%%%%%%%%%%%%%%%%%%%%
\subsubsection{Triangular Grid}
\label{eb_njk_typical_tria}
%%%%%%%%%%%%%%%%%%%%%%%%%%%%%%%%%%%%%%%%%%%%%%%%

 \indent 
 
\begin{figure}[htbp!]
    \centering
          \begin{subfigure}[t]{0.48\textwidth}
    \centering
        \includegraphics[width=0.7\textwidth,trim=0 0 0 0,clip]{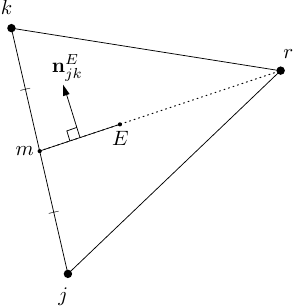}
          \caption{Triangular element.}
          \label{fig:njk_contribution_2d}
      \end{subfigure}
      \hfill
          \begin{subfigure}[t]{0.48\textwidth}
    \centering
        \includegraphics[width=0.7\textwidth,trim=0 0 0 0,clip]{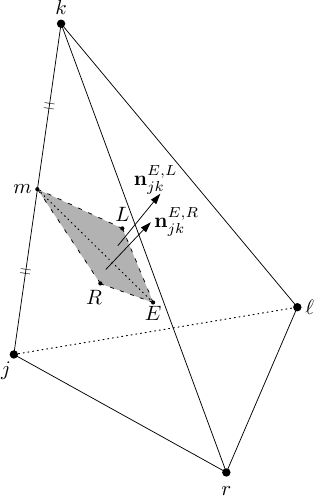}
          \caption{Tetrahedral element.}
          \label{fig:njk_contribution_3d}
      \end{subfigure}
            \caption{
\label{fig:njk_contribution}%
Lumped directed-area vector contributions from an element $E$ to the edge $\{ j,k \}$.} 
\end{figure}

Typically, the lumped directed-area vectors are computed by summing dual-face contributions at edges in a loop over elements. For a triangular grid, at an edge $\{ j,k \}$ in an element $E$ as shown in Figure \ref{fig:njk_contribution_2d}, we form a dual face by connecting the edge midpoint $(x_m,y_m) = \left(  \frac{ x_j + x_k }{2}, \frac{ y_j + y_k  }{2} \right)$ and the geometric centroid of the element $(x_E,y_E) =   \left(  \frac{ x_j + x_k + x_r}{3}, \frac{ y_j + y_k  + y_r}{3} \right)$, and compute the directed-area vector: 
\begin{eqnarray}
%   {\bf n}_{jk}^E = (  y_m  - y_E  , x_E - x_m    ), 
   {\bf n}_{jk}^E = (  \Delta y_{mE}  ,  \Delta x_{Em}     ), 
\end{eqnarray}
where $ \Delta y_{mE} = y_m  - y_E $, $\Delta x_{Em} = x_E - x_m$, and we have assumed the edge is oriented from $j$ to $k$ (i.e., $j < k$). Then, we add it to the lumped directed-area vector of the edge $\{ j,k \}$: 
\begin{eqnarray}
 {\bf n}_{jk} =  {\bf n}_{jk}  +  {\bf n}_{jk}^E.
\end{eqnarray}
Repeat the same for the other two edges, $\{ j,r \}$ and $\{ r,k \}$, and move onto the next element. At the end of the process, we will have the lumped directed-area vectors at all edges. 

%%%%%%%%%%%%%%%%%%%%%%%%%%%%%%%%%%%%%%%%%
%  
%%%%%%%%%%%%%%%%%%%%%%%%%%%%%%%%%%%%%%%%%
%%%%%%%%%%%%%%%%%%%%%%%%%%%%%%%%%%%%%%%%%
% 
%%%%%%%%%%%%%%%%%%%%%%%%%%%%%%%%%%%%%%%%%

%%%%%%%%%%%%%%%%%%%%%%%%%%%%%%%%%%%%%%%%%%%%%%%%
\subsubsection{Tetrahedral Grid}
\label{eb_njk_typical_tet}
%%%%%%%%%%%%%%%%%%%%%%%%%%%%%%%%%%%%%%%%%%%%%%%%

\indent

Similarly, for a tetrahedral grid, we typically compute the lumped directed-area vectors by summing dual-face contributions at edges in a loop over tetrahedra. In this case, we have two triangular dual faces in each element: $\{  {\bf x}_m,  {\bf x}_E,  {\bf x}_L \}$ and $\{  {\bf x}_m, {\bf x}_R,  {\bf x}_E \}$ as illustrated in Figure \ref{fig:njk_contribution_3d}. Here, $m$ is the midpoint of the edge $\{ j,k \}$, $E$ is the geometric centroid of the element $E$, and $L$ and $R$ are the geometric centroids of the triangular faces $\{ j,\ell,k \}$ and $\{ j,r,k \}$, respectively: 
%: $(x_E,y_E,z_E) =   \left(  \frac{ x_j + x_k + x_\ell + x_r}{4}, \frac{ y_j + y_k + y_\ell + y_r}{4}, \frac{ z_j + z_k + z_\ell + z_r}{4}  \right)$, $L$ and $R$ are the geometric centroids of the triangular faces $[j,\ell,k]$ and $[j,r,k]$, respectively: $(x_L, y_L, z_L) =   \left(  \frac{ x_j + x_\ell + x_k}{3}, \frac{ y_j + y_\ell  + y_k}{3} \right)$ and 
%$(x_R, y_R, z_R) =   \left(  \frac{ x_j + x_r + x_k}{3}, \frac{ y_j + y_r  + y_k}{3}, \frac{ z_j + z_r  + z_k}{3} \right)$. 
\begin{eqnarray} 
(x_m,y_m,z_m)  &=&  \left(  \frac{ x_j + x_k }{2}, \frac{ y_j + y_k  }{2}, \frac{ z_j + z_k  }{2} \right), \\ [2ex]
(x_E,y_E,z_E)   &=&   \left(  \frac{ x_j + x_k + x_\ell + x_r}{4}, \frac{ y_j + y_k + y_\ell + y_r}{4}, \frac{ z_j + z_k + z_\ell + z_r}{4}  \right), \\ [2ex]
(x_L, y_L, z_L)  &=&  \left(  \frac{ x_j + x_\ell + x_k}{3}, \frac{ y_j + y_\ell  + y_k}{3}, \frac{ z_j + z_\ell  + z_k}{3}  \right), \\ [2ex]
(x_R, y_R, z_R)  &=&   \left(  \frac{ x_j + x_r + x_k}{3}, \frac{ y_j + y_r  + y_k}{3}, \frac{ z_j + z_r  + z_k}{3} \right).
\end{eqnarray}
Therefore, we add the directed-area vectors of the two triangular dual faces ${\bf n}_{jk}^{E,L}$ and ${\bf n}_{jk}^{E,R}$ to the lumped directed-area vector ${\bf n}_{jk}$: 
\begin{eqnarray}
% {\bf n}_{jk} =  {\bf n}_{jk}  +  \left( {\bf n}_{jk}^{E,L} + {\bf n}_{jk}^{E,R} \right),
 {\bf n}_{jk} =  {\bf n}_{jk}  +  {\bf n}_{jk}^{E,L} + {\bf n}_{jk}^{E,R} ,
\end{eqnarray}
where
\begin{eqnarray}
  {\bf n}_{jk}^{E,L}  &=&  \frac{1}{2} \left(   \Delta y_{Em} \Delta z_{Lm}  - \Delta z_{Em}  \Delta y_{Lm},
                                                                  \Delta z_{Em} \Delta x_{Lm}  - \Delta x_{Em}  \Delta z_{Lm},
                                                                  \Delta x_{Em} \Delta y_{Lm}  - \Delta y_{Em}  \Delta x_{Lm}    \right) ,  \\ [2ex]
  {\bf n}_{jk}^{E,R}  &=&  \frac{1}{2}  \left(   \Delta y_{Rm} \Delta z_{Em}  - \Delta z_{Rm}  \Delta y_{Em},
                                                                  \Delta z_{Rm} \Delta x_{Em}  - \Delta x_{Rm}  \Delta z_{Em},
                                                                  \Delta x_{Rm} \Delta y_{Em}  - \Delta y_{Rm}  \Delta x_{Em}  \right) ,
\end{eqnarray}
where again we have assumed $j < k$, so that the edge is oriented from $j$ to $k$. We repeat the same for the other five edges in the element $E$, and move onto the next element. At the end of the process, we will have the lumped directed-area vectors at all edges. In the actual implementation, the factor $1/2$ can be applied, after summing $2 {\bf n}_{jk}^{E,L}$ and $2 {\bf n}_{jk}^{E,R}$ over elements, in a subsequent loop over edges for efficiency.

% Elements
% 3 additions and 3 multiplications for (xc,yc,zc)
% 6 additions and 3 multiplications for face centroids, 4 faces ->  24 add + 12 mul
%------
% 27 add + 15 mul
% Edges (6)
%3 additions and 3 multiplications for (xm,ym,zm)
%6 additions and 6 multiplications for computing ${\bf n}_{jk}^{E,L}$ (the factor $1/2$ ignored here), 
%6 additions and 6 multiplications for computing ${\bf n}_{jk}^{E,R}$ (the factor $1/2$ ignored here), 
%3 additions for adding them, and 3 additions for adding the contributions to ${\bf n}_{jk}$. 
%-------
% 21 add + 15 mul ->  x 6 edges = 126 add + 90 mul
%
% 153 add        + 105 mul  per element
% 153*NE add  + 105*NE mul in total

%%%%%%%%%%%%%%%%%%%%%%%%%%%%%%%%%%%%%%%%%%%%%%%%%%
\begin{figure}[h!]
\begin{center}
\begin{minipage}[b]{0.95\textwidth}
\begin{center}
        \includegraphics[width=0.5\textwidth,trim=0 0 0 0,clip]{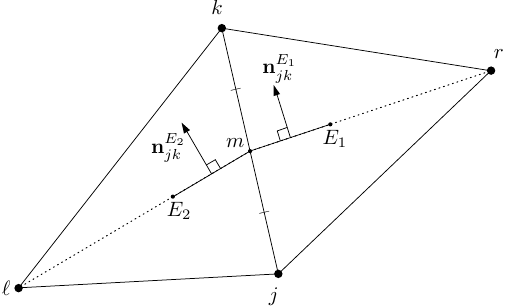}
          \caption{Lumped directed-area vector contributions from two triangles to the edge $\{ j,k \}$.}
          \label{fig:efficient_njk_2d}
\end{center}
\end{minipage}
\end{center}
\end{figure}
%%%%%%%%%%%%%%%%%%%%%%%%%%%%%%%%%%%%%%%%%%%%%%%%%%

%%%%%%%%%%%%%%%%%%%%%%%%%%%%%%%%%%%%  
%%%%%%%%%%%%%%%%%%%%%%%%%%%%%%%%%%%%  
%%%%%%%%%%%%%%%%%%%%%%%%%%%%%%%%%%%%  
\subsection{An Efficient Algorithm for Computing Lumped Directed-Area Vectors}
\label{eb_njk_vj_typical}
%%%%%%%%%%%%%%%%%%%%%%%%%%%%%%%%%%%%  
%%%%%%%%%%%%%%%%%%%%%%%%%%%%%%%%%%%%  
%%%%%%%%%%%%%%%%%%%%%%%%%%%%%%%%%%%% 

%%%%%%%%%%%%%%%%%%%%%%%%%%%%%%%%%%%%%%%%%%%%%%%%
\subsubsection{Triangular Grid}
\label{eb_njk_efficient_tria}
%%%%%%%%%%%%%%%%%%%%%%%%%%%%%%%%%%%%%%%%%%%%%%%%

 \indent 
 
For a triangular grid, we derive a more efficient algorithm for computing the lumped directed-area vectors as follows. Consider the triangle $E_1$ having the edge $\{ j,k \}$ as illustrated in Figure \ref{fig:efficient_njk_2d}. By summing the outward directed-area vectors of the three faces of the triangle $\{ m,r,k \}$, we have
\begin{eqnarray}
3 ( - {\bf n}_{jk}^{E_1} )  +      {\bf n}_j^{E_1}  +  \frac{1}{2} {\bf n}_r  =  {\bf 0},  
%    {\bf n}_j   +   \frac{1}{2}  {\bf n}_2  +  \frac{1}{2} {\bf n}_3  =  {\bf 0},   
\end{eqnarray}
where ${\bf n}_j^{E_1} = (   y_k - y_r   , x_r-x_k  )$ is the outward directed-area vector of the face opposite to the node $j$, ${\bf n}_r = (  y_j - y_k   , x_k - x_j   )$ is the outward directed-area vector of the face opposite to the node $r$. It gives
\begin{eqnarray}
 {\bf n}_{jk}^{E_1}  = \frac{1}{3}  \left[  {\bf n}_j^{E_1} +  \frac{1}{2} {\bf n}_r  \right].
 \label{njk_E1_2d}
\end{eqnarray}
Similarly, for the triangle $E_2$, we obtain
\begin{eqnarray}
 {\bf n}_{jk}^{E_2}  = \frac{1}{3}  \left[  {\bf n}_j^{E_2}   +  \frac{1}{2} {\bf n}_\ell  \right].
 \label{njk_E2_2d}
\end{eqnarray}
The lumped directed-area vector for the edge $\{ j,k \}$ is then given by the sum of the two: 
 \begin{eqnarray}
  {\bf n}_{jk} = \frac{1}{3} (   {\bf n}_{j}^{E_1}   +    {\bf n}_{j}^{E_2} ),
\end{eqnarray}
where ${\bf n}_r$ and ${\bf n}_\ell$ have been canceled out because ${\bf n}_\ell  = - {\bf n}_r $. It can be formally expressed as
 \begin{eqnarray}
  {\bf n}_{jk} = \frac{1}{3} \sum_{ E \in \{ E_{jk} \} }   {\bf n}_{j}^{E} ,
  \label{njk_eff_tria}
\end{eqnarray}
where $\{ E_{jk} \} $ denotes the set of elements sharing the edge $\{ j,k \}$. %Note that this formula does not require the construction of dual faces and only need a directed-area vector of a face opposite to a node.

The formula (\ref{njk_eff_tria}) can be implemented in a loop over elements: distribute the local contribution as 
 \begin{eqnarray}
  {\bf n}_{jk} = {\bf n}_{jk}  +   \frac{1}{3}   {\bf n}_{j}^{E} ,
  \label{njk_eff_algorithm_tria_01}
\end{eqnarray}
if the edge $\{ j,k \}$ is oriented from $j$ to $k$ (i.e., if $j < k$), or 
 \begin{eqnarray}
  {\bf n}_{kj} = {\bf n}_{kj}  +   \frac{1}{3}   {\bf n}_{k}^{E} ,
  \label{nkj_eff_algorithm_tria_02}
\end{eqnarray}
if the edge $\{ j,k \}$ is oriented from $k$ to $j$ (i.e., if $j > k$), perform the same for all three edges within an element, and move onto the next element. 

\begin{figure}[t]
    \centering 
          \begin{subfigure}[t]{0.48\textwidth}
    \centering
        \includegraphics[width=1.0\textwidth,trim=0 0 0 0,clip]{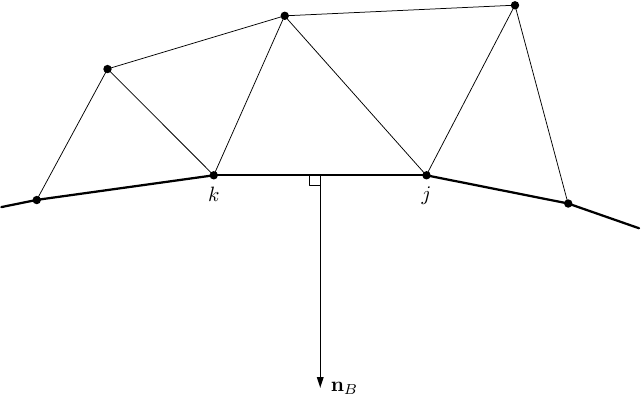}
          \caption{Boundary edge visited in a boundary edge loop.}
          \label{fig:njk_contribution_2d_b_collapsed_bgrid}
      \end{subfigure}
      \hfill
          \begin{subfigure}[t]{0.48\textwidth}
    \centering
        \includegraphics[width=1.0\textwidth,trim=0 0 0 0,clip]{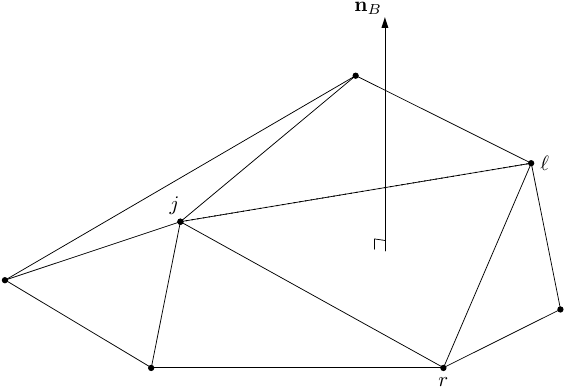}
          \caption{Boundary triangle $B$  visited in a boundary element loop.}
          \label{fig:njk_contribution_3d_b_loop}
      \end{subfigure}
            \caption{
 \label{fig:njk_contribution_3d_b}
Boundary elements in two and three dimensions. Note that the directed-area vector ${\bf n}_B$ of a boundary element is pointing outward from the interior of a domain in both figures. } 
\end{figure}

If an edge is a boundary edge, e.g., the edge $\{ j,k \}$ has only $E_1$ in Figure \ref{fig:efficient_njk_2d}, then the term $ \frac{1}{2} {\bf n}_r$ in Equation (\ref{njk_E1_2d}) needs to be added. This boundary correction can be implemented easily by making a loop over boundary edges, $B=1,2,3, \cdots, N_B$, where $N_B$ is the number of boundary edges and performing the following correction for each boundary edge $B$ with two ends nodes $j$ and $k$ (see Figure \ref{fig:njk_contribution_2d_b_collapsed_bgrid}), 
 \begin{eqnarray}
  {\bf n}_{jk} =   {\bf n}_{jk} + \frac{1}{6} {\bf n}_{B}, 
  \label{njk_eff_algorithm_tria_02p5_nB}
\end{eqnarray}
where ${\bf n}_{B}$ is the outward directed-area vector of the boundary edge $\{ j,k \}$. Note that this boundary correction is independent of the orientation of the edge $\{ j , k \}$, i.e., we add the same $\frac{1}{6} {\bf n}_{B}$: 
 \begin{eqnarray}
  {\bf n}_{kj} =   {\bf n}_{kj} + \frac{1}{6} {\bf n}_{B}, 
  \label{njk_eff_algorithm_tria_02p5_nB_21}
\end{eqnarray}
even if the edge $\{ j, k \}$ is actually oriented from $k$ to $j$. To see this, we consider Figure \ref{fig:njk_contribution_2d} and assume the edge $\{ j, k \}$ is a boundary edge. Summing over the outward directed-area vectors of the three faces of the triangle $\{ m,r,k \}$ and also of the triangle $\{ m, r, j \}$, we obtain
\begin{eqnarray}
3 ( - {\bf n}_{jk}^{E} )  +      {\bf n}_j^{E}  +  \frac{1}{2} {\bf n}_r   &=&  {\bf 0},  \\ [2ex]
3 ( -  {\bf n}_{kj}^{E} )  +      {\bf n}_k^{E}  +  \frac{1}{2} {\bf n}_r  &=&  {\bf 0},  
%    {\bf n}_j   +   \frac{1}{2}  {\bf n}_2  +  \frac{1}{2} {\bf n}_3  =  {\bf 0},   
\end{eqnarray}
respectively, where ${\bf n}_{kj}^{E} = - {\bf n}_{jk}^{E}$, and find
\begin{eqnarray}
   {\bf n}_{jk}^{E}   &=&  \frac{1}{3}  \left[  {\bf n}_j^{E}  +  \frac{1}{2} {\bf n}_r  \right],  \\ [2ex]
   {\bf n}_{kj}^{E}   &=&  \frac{1}{3}  \left[  {\bf n}_k^{E} +  \frac{1}{2} {\bf n}_r  \right].
\end{eqnarray}
Now turn to Figure \ref{fig:njk_contribution_2d_b_collapsed_bgrid}, and find ${\bf n}_r = {\bf n}_B$ and $ {\bf n}_{jk} = {\bf n}_{jk}^{E}$ or ${\bf n}_{kj} = {\bf n}_{kj}^{E}$ (there is only one element for the edge). If the edge $\{ j, k \}$ is oriented from $j$ to $k$, we have ${\bf n}_{jk}$, the first term is taken into account in Equation (\ref{njk_eff_algorithm_tria_01}), and the required correction is $\frac{1}{6} {\bf n}_B$. On the other hand, if the edge $\{ j, k \}$ is oriented from $k$ to $j$, we have ${\bf n}_{kj}$, the first term is taken into account in Equation (\ref{nkj_eff_algorithm_tria_02}), and the required correction is the same $\frac{1}{6} {\bf n}_B$.
 
% $22 N_E$ additions and $8 N_E$ multiplications,

% n2 = 3*n12 + nR 
% Note that (-n1) + (-n2) + 2*nR = 0 ->  n2 = -n1 + 2*nR,   
% So, (-n1) + 2*nR = 3*n12 + nR
% n12 = 1/3 * [ (-n1) +  nR ]

% 147 N_E + 3 N_e ~ 147 N_E + 3*(7/6*NE) =  150.5*NE
% 36  N_E + 3 N_e ~ 36 NE + 3*(7/6*NE)      =    39.5*NE
%

%  
%%%%%%%%%%%%%%%%%%%%%%%%%%%%%%%%%%%%%%%%%
%  \begin{figure}[htbp!]
%%%%%%%%%%%%%%%%%%%%%%%%%%%%%%%%%%%%%%%%%
% 
%%%%%%%%%%%%%%%%%%%%%%%%%%%%%%%%%%%%%%%%%

%%%%%%%%%%%%%%%%%%%%%%%%%%%%%%%%%%%%%%%%%%%%%%%%
\subsubsection{Tetrahedral Grid}
\label{eb_njk_efficient_tet}
%%%%%%%%%%%%%%%%%%%%%%%%%%%%%%%%%%%%%%%%%%%%%%%%

\indent

For a tetrahedral grid, we consider local contributions to a node $j$ within an element $E$. See Figure \ref{fig:njk_contribution_3d}. 
By summing the outward directed-area vectors of the tetrahedron $\{ m,r,\ell,k \}$, we have
\begin{eqnarray}
% 3/2*nR + 3/2*nL + (-n4) + 6*n14 = 0  -- (1)
%  \frac{1}{2}  {\bf n}_2 +  \frac{1}{2}  {\bf n}_3 + {\bf n}_k +  6 {\bf n}_{jk}^{E}  = {\bf 0},  \\ [2ex]
% 3*(nL + nR) + (-n4) + (-n1) 
%   {\bf n}_2 +  {\bf n}_3 + {\bf n}_k +    {\bf n}_j   = {\bf 0}.
6 ( - {\bf n}_{jk}^{E,L} - {\bf n}_{jk}^{E,R}  )  + 
    {\bf n}_j^E   +   \frac{1}{2}  {\bf n}_\ell  +  \frac{1}{2} {\bf n}_r  =  {\bf 0},   
\end{eqnarray}
where ${\bf n}_j^E$ is the outward directed-area vector of the triangular face opposite to the node $j$ in the element $E$, and similarly for $ {\bf n}_\ell^E$ and $ {\bf n}_r^E$. It gives 
\begin{eqnarray} 
  %\equiv
   {\bf n}_{jk}^{E,L} + {\bf n}_{jk}^{E,R} 
  = \frac{1}{6}  
\left[
    {\bf n}_j^E   +   \frac{1}{2} (  {\bf n}_\ell^E  +  {\bf n}_r^E ) 
\right],
\label{njk_local_formula_3d}
\end{eqnarray} 
which is the total local contribution to the lumped directed-area vector ${\bf n}_{jk}$ (the edge $\{ j,k \}$ is assumed to be oriented from $j$ to $k$). Note that the dual face $\{ m,R,E,L \}$ is planar, the geometric centroid $E$ is located halfway between $m$ and the midpoint of the edge $\{ r,\ell \}$, and it leads to $6 ( - {\bf n}_{jk}^{E,L} - {\bf n}_{jk}^{E,R}  )$ as the outward directed-area vector of the triangle $\{ m,r,\ell \}$. Summing the local contributions from all tetrahedra sharing the edge $\{ j,k \}$, we obtain
 \begin{eqnarray}
  {\bf n}_{jk} 
  = \sum_{ E \in \{ E_{jk} \} }  \frac{1}{6}  
\left[
    {\bf n}_j^E   +   \frac{1}{2} (  {\bf n}_\ell^E  +  {\bf n}_r^E ) 
    \right]
  = \frac{1}{6} \sum_{ E \in \{ E_{jk} \} }   {\bf n}_{j}^{E} ,
  \label{efficient_njk_3d}
\end{eqnarray}
where $\{ E_{jk} \} $ denotes the set of elements sharing the edge $\{ j,k \}$. Note that we consider interior edges here (i.e., there is always an adjacent tetrahedron to the face opposite to the node $\ell$ and the same for the node $r$), and therefore, the term $\frac{1}{2} (  {\bf n}_\ell^E  +  {\bf n}_r^E ) $ has canceled out. 

As in the triangular-grid case, we implement this formula in an element loop, adding $ 2 {\bf n}_{j}^{E}$ to $ {\bf n}_{jk}$ if the edge is oriented from $j$ to $k$ (i.e., if $j < k$): 
 \begin{eqnarray}
  {\bf n}_{jk} 
=
    {\bf n}_{jk} 
  + 
%  \left(   \Delta y_{k \ell} \Delta z_{k \ell}  - \Delta z_{k \ell}  \Delta y_{k \ell},
%            \Delta z_{k \ell} \Delta x_{k \ell}  - \Delta x_{k \ell}  \Delta z_{k \ell},
%            \Delta x_{k \ell} \Delta y_{k \ell}  - \Delta y_{k \ell}  \Delta x_{k \ell}    \right),
  \left(   \Delta y_{r k } \Delta z_{\ell k  }  - \Delta z_{r k  }  \Delta y_{\ell k  },
            \Delta z_{r k } \Delta x_{\ell k  }  - \Delta x_{r k  }  \Delta z_{\ell k  },
            \Delta x_{r k } \Delta y_{\ell k  }  - \Delta y_{r k  }  \Delta x_{\ell k  }    \right),
  \label{add_contribution_to_njk_tetra}
\end{eqnarray}
or adding $ 2 {\bf n}_{k}^{E}$ to $ {\bf n}_{kj}$ if the edge is oriented from $k$ to $j$ (i.e., if $j > k$): 
 \begin{eqnarray}
  {\bf n}_{kj} 
=
    {\bf n}_{kj} 
  + 
%  \left(   \Delta y_{k \ell} \Delta z_{k \ell}  - \Delta z_{k \ell}  \Delta y_{k \ell},
%            \Delta z_{k \ell} \Delta x_{k \ell}  - \Delta x_{k \ell}  \Delta z_{k \ell},
%            \Delta x_{k \ell} \Delta y_{k \ell}  - \Delta y_{k \ell}  \Delta x_{k \ell}    \right),
  \left(   \Delta y_{r j } \Delta z_{\ell j  }  - \Delta z_{r j  }  \Delta y_{\ell j  },
            \Delta z_{r j } \Delta x_{\ell j  }  - \Delta x_{r j  }  \Delta z_{\ell j  },
            \Delta x_{r j } \Delta y_{\ell j  }  - \Delta y_{r j  }  \Delta x_{\ell j  }    \right),
  \label{add_contribution_to_nkj_tetra}
\end{eqnarray}
and doing the same for the other five edges in an element $E$, followed by a subsequent loop over edges, where we apply the factor $1/12$:
 \begin{eqnarray}
    {\bf n}_{jk} 
= 
  \frac{1}{12} \,  {\bf n}_{jk} ,
  \label{apply_1over12_at_the_end_njk_tetra}
\end{eqnarray}
where the factor $1/12$ is the product of $1/6$ in the formula (\ref{efficient_njk_3d}) and $1/2$ to cancel $2$ of $2 {\bf n}_{j}^{E}$ or $2 {\bf n}_{k}^{E}$. 

If the edge $\{ j,k \}$ is a boundary edge, the term $\frac{1}{2} (  {\bf n}_\ell^E  +  {\bf n}_r^E ) $ does not cancel out and needs to be added. Similarly to the triangular-grid case, we can add this contribution by looping over boundary elements, $B=1,2,3, \cdots, N_B$, where $N_B$ is the number of boundary elements, and performing the following correction for each boundary triangle $B$ with nodes $j$, $r$, $\ell$ (see Figure \ref{fig:njk_contribution_2d_b_collapsed_bgrid}), 
 \begin{eqnarray}
    {\bf n}_{j r} =  {\bf n}_{j r} + {\bf n}_B, \quad
    {\bf n}_{r \ell} =  {\bf n}_{r \ell } + {\bf n}_B, \quad
    {\bf n}_{\ell j} =  {\bf n}_{\ell j} + {\bf n}_B, 
    \label{njk_boundary_correction_3d}
\end{eqnarray}
where $ {\bf n}_{j r}$,  $ {\bf n}_{r \ell}$,  and $ {\bf n}_{\ell j}$ are the lumped directed-area vectors of the edges $\{ j, r \}$, $\{ k, \ell \}$, and $\{ \ell, j \}$, respectively, and ${\bf n}_B$ is the outward directed-area of the boundary face $B$:
 \begin{eqnarray}
  {\bf n}_{B} 
=
\frac{1}{2}
   \left(  \Delta y_{r j } \Delta z_{\ell j  }  - \Delta z_{r j  }  \Delta y_{\ell j  },
            \Delta z_{r j} \Delta x_{ \ell j  }  - \Delta x_{r j  }  \Delta z_{\ell j  },
            \Delta x_{r j } \Delta y_{\ell j  }  - \Delta y_{r j  }  \Delta x_{\ell j  }    \right).
  \label{tetra_tria_boundary_nB}
\end{eqnarray}
 %The factor $1/2$ missing here is then applied in the next edge loop for Equation (\ref{apply_1over12_at_the_end_njk_tetra}). 
Note, as in the two-dimensional case, that the boundary correction (\ref{njk_boundary_correction_3d}) is independent of the edge orientation, and thus ${\bf n}_{jr}$ may actually be $ {\bf n}_{rj}$ directed from $r$ to $j$, for example.   

% Ne = 7*Nv  , NE = 6*Nv

%2D:
%  Nv - Ne + NE = 2
%  Nv - Nei - Neb + NE = 2
%    2*Nei + Neb = 3*NE  (tria)
%  Nv - 3/2*NE + Neb/2 - Neb + NE = 2
%  Nv - 1/2*NE - Neb/2 = 2
%   -->  NE = 2*Nv - Neb - 4
%   Ne = NE + Nv - 2 = 3*Nv - Neb - 6

%3D: 
%  Nv - Ne + Nf - NE = 0
% Nv - Ne + Nfi+Nfb - NE = 0
%    2*Nfi + Nfb = 4*NE  (tetra)
% Nv - Ne + 2*NE - Nfb/2 +Nfb - NE = 0
% Nv - Ne + NE + Nfb/2 = 0
% NE + Nv + Nfb/2 = Ne 
% If NE = 6*Nv, then
%    Ne = 7*Nv + Nfb/2

% n12 + n13 + n14 + nL + nR + nB = 0
% 3*(nL + nR + nB) + (-n1) = 0 

% S - 2/3*1/2*S = nR   -> 2/3*S = nR -> S = 3/2*nR. 
% 3/2*nR + 3/2*nL + (-n4) + 6*n14 = 0  -- (1)
% 3*(nL + nR) + (-n4) + (-n1) = 0   --(2)
% 3/2*nR + 3/2*nL + n1 - 3(nL+nR) + 6*n14 = 0
% 6*n14 = (-n1) + 3/2*(nL+nR)
% n14 = 1/6 * (  (-n1) + 3/2*(nL + nR )  )
% n14 = 1/6 * (  (-n1) + 3/2*(nL + nR )  )

% 3*(nL + nR) + (-n4) + (-n1) = 0 
%      We have in AIAA2010-5093 n4 = 6*n14 + ( n1 + n4 )/2 
% ( n4 - n1 )/2 = 6 *n14

%%%%%%%%%%%%%%%%%%%%%%%%%%%%%%%%%%%%%%%%%
%  
%%%%%%%%%%%%%%%%%%%%%%%%%%%%%%%%%%%%%%%%%
 % \begin{figure}[htbp!]
  \begin{figure}[t]
    \centering
          \begin{subfigure}[t]{0.48\textwidth}
    \centering
        \includegraphics[width=0.5\textwidth,trim=0 0 0 0,clip]{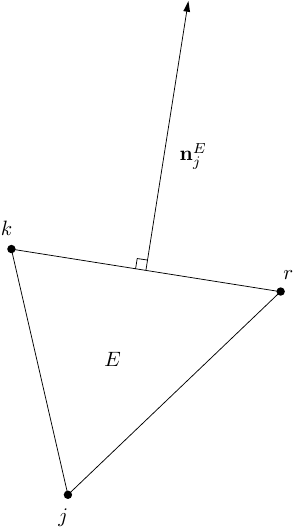}
          \caption{Triangular element.}
          \label{fig:njk_contribution_2d_njE}
      \end{subfigure}
      \hfill
          \begin{subfigure}[t]{0.48\textwidth}
    \centering
        \includegraphics[width=0.6\textwidth,trim=0 0 0 0,clip]{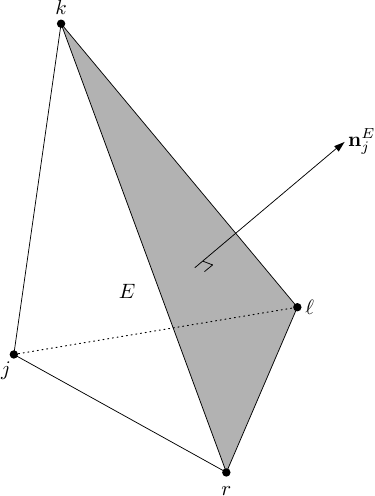}
          \caption{Tetrahedral element.}
          \label{fig:njk_contribution_3d_njE}
      \end{subfigure}
            \caption{
\label{fig:njk_contribution_njE}%
Directed-area vectors needed to compute ${\bf n}_{jk}$ oriented from $j$ to $k$. } 
\end{figure}
%%%%%%%%%%%%%%%%%%%%%%%%%%%%%%%%%%%%%%%%%
% 
%%%%%%%%%%%%%%%%%%%%%%%%%%%%%%%%%%%%%%%%%

%%%%%%%%%%%%%%%%%%%%%%%%%%%%%%%%%%%%%%%%%%%%%%%%
\subsubsection{General Form and Boundary Closure}
\label{eb_njk_general}
%%%%%%%%%%%%%%%%%%%%%%%%%%%%%%%%%%%%%%%%%%%%%%%%

\indent

The two- and three-dimensional formulas can be written in the following general form: 
 \begin{eqnarray}
  {\bf n}_{jk} = \frac{2}{D ( D + 1) } \sum_{ E \in \{ E_{jk} \} }   {\bf n}_{j}^{E},
  \label{efficient_njk_formula_general_D}
\end{eqnarray}
where $D=2$ for a triangular grid and $D=3$ for a tetrahedral grid, with a boundary correction:
 \begin{eqnarray}
    {\bf n}_{jk} =  {\bf n}_{jk} +  \frac{2}{D ( D + 1) }  \frac{1}{2} {\bf n}_B = {\bf n}_{jk} +  \frac{1}{D ( D + 1) }  {\bf n}_B ,  
\end{eqnarray}
for all edges of a boundary element $B$ (1 edge in a triangular grid and 3 edges in a tetrahedral grid), for $B = 1,2,3, \cdots, N_B$, where ${\bf n}_B$ is the outward directed-area vector of a boundary element $B$. Here is an important observation: these formulas allow us to compute the lumped directed-area vector ${\bf n}_{jk}$ without actually forming the dual faces: it only requires the directed-area vector of the face opposite to the node $j$ for each element $E \in \{ E_{jk} \}$, as illustrated in Figure \ref{fig:njk_contribution_njE}, and outward directed-area vectors of each boundary element. Therefore, there is no need at all to form the median dual faces and control volumes (which requires edge midpoints, face centroids, element centroids) in the implementation of the edge-based discretization. This is the main reason that the proposed algorithm is more efficient and significantly simpler to implement than the traditional algorithm, especially for tetrahedral grids. 

\begin{figure}[t]
    \centering
          \begin{subfigure}[t]{0.48\textwidth}
    \centering
        \includegraphics[width=0.8\textwidth,trim=0 0 0 0,clip]{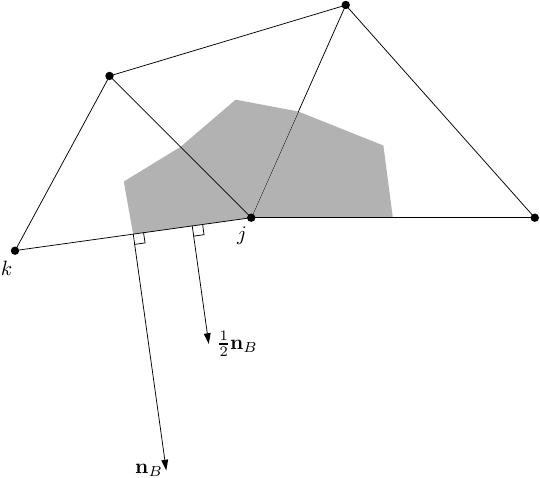}
          \caption{Boundary of a triangular grid.}
          \label{fig:boundary_closure_2d}
      \end{subfigure}
      \hfill
          \begin{subfigure}[t]{0.48\textwidth}
    \centering
        \includegraphics[width=0.8\textwidth,trim=0 0 0 0,clip]{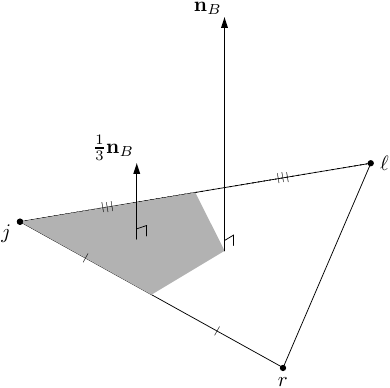}
          \caption{Boundary element of a tetrahedral grid.}
          \label{fig:boundary_closure_3d}
      \end{subfigure}
            \caption{
          \label{fig:boundary_closure}%
Outward directed-area vectors at boundary elements in two and three dimensions. In both cases, ${\bf n}_B$ is pointing outward from the interior domain. The shaded area is the median dual volume around the node $j$ in the triangular grid (left); the shaded area is a boundary face of the median dual volume around the node $j$ in the tetrahedral grid case (right).} 
\end{figure}

It is always good to verify the implementation of an algorithm in a code. One simple verification for the lumped directed-area vectors is to confirm that the sum of the lumped directed-area vectors is zero at every node, which can be performed as follows: set ${\bf a}_j  = {\bf 0}$ at all nodes ($ j = 1,2,3, \cdots , N_v$), and accumulate the lumped directed-area vectors in a loop over edges, 
 \begin{eqnarray}
 {\bf a}_j = {\bf a}_j + {\bf n}_{jk}, \quad 
 {\bf a}_k = {\bf a}_k - {\bf n}_{jk}, 
\end{eqnarray}
where the edge $\{ j,k \}$ is assumed to be oriented from $j$ to $k$. At the end of the process, we must have ${\bf a}_j = {\bf 0}$ at all interior nodes. At boundary nodes, ${\bf a}_j \ne {\bf 0}$ at this point, but we must have zero sums at boundary nodes as well in order to correctly discretize a conservation law at boundary nodes: the flux balance must be zero if the flux is constant. The edge-based discretization (\ref{threed_fv_semidiscrete_system_00}) actually needs to be closed at boundary elements with an accuracy-preserving boundary flux quadrature formula as discussed in Ref.~\cite{nishikawa:AIAA2010} for second-order accuracy and Ref.~\cite{nishikawa_boundary_quadrature:JCP2015} for third-order accuracy. For verification of ${\bf n}_{jk}$, it suffices to add outward directed-area vectors corresponding to dual faces at boundary elements, each of which is precisely $1/D$ of the outward directed-area vector of a boundary element: $1/2 \, {\bf n}_B = 1/2 ( y_j - y_k, x_k - x_j)$ as illustrated in Figure \ref{fig:boundary_closure_2d} and 
$1/3 \, {\bf n}_B = 1/6 (   
             \Delta y_{r j } \Delta z_{\ell j}  - \Delta z_{r j  }  \Delta y_{\ell  j  },
            \Delta z_{ r j } \Delta x_{\ell j  }  - \Delta x_{r j  }  \Delta z_{\ell  j  },
            \Delta x_{ r j  } \Delta y_{\ell j   }  - \Delta y_{r j  }  \Delta x_{\ell j  } )$ as illustrated in Figure \ref{fig:boundary_closure_3d}. 
Therefore, for each node $i_B$ of a boundary element $B$, $B=1,2,3, \cdots, N_B$, we perform the following:
 \begin{eqnarray}
 {\bf a}_{i_B}= {\bf a}_{i_B} +  \frac{1}{D} {\bf n}_B,  
\end{eqnarray}
where $i_B = j $ and $k$ in two dimensions (see Figure \ref{fig:boundary_closure_2d}) and $ i_B = j$, $r$, and $\ell$ in three dimensions (see Figure \ref{fig:boundary_closure_3d}). At the end of this loop, we must have ${\bf a}_j = {\bf 0}$ at all boundary nodes. 

%%%%%%%%%%%%%%%%%%%%%%%%%%%%%%%%%%%%%%%%%
%  
%%%%%%%%%%%%%%%%%%%%%%%%%%%%%%%%%%%%%%%%%
 % \begin{figure}[htbp!]
%%%%%%%%%%%%%%%%%%%%%%%%%%%%%%%%%%%%%%%%%
% 
%%%%%%%%%%%%%%%%%%%%%%%%%%%%%%%%%%%%%%%%%

%But numerical results indicate that the 

%%%%%%%%%%%%%%%%%%%%%%%%%%%%%%%%%%%% 
\subsection{Complexity and Discussion}
\label{eb_njk_complexity}
%%%%%%%%%%%%%%%%%%%%%%%%%%%%%%%%%%%% 

It is not simple to precisely measure the relative efficiency of the two algorithms as it depends on various factors, e.g., memory access, compilers, code optimizations, and architectures. However, it is often helpful to provide a complexity analysis based on the number of additions and multiplications required in each algorithm and compare them, ignoring all other factors (subtractions are taken as additions, and divisions of constant values are assumed to be implemented as multiplications). For simplicity, we ignore boundary effects (i.e., assuming large grids) and focus on interior edges and nodes. Later, we will present comparison of actual computing times and show that the actual speed-up factors are not too far from those indicated by the complexity analysis.

%%%%%%%%%%%%%%%%%%%%%%%%%%%%%%%%%%%% 
\subsubsection{Triangular Grid}
\label{eb_njk_complexity_tria}
%%%%%%%%%%%%%%%%%%%%%%%%%%%%%%%%%%%% 

For a triangular grid, we first consider the traditional algorithm described in Section \ref{eb_njk_typical_tria}. Here, we first need 2 additions and 2 multiplications for computing the edge midpoint coordinates. Then, it takes 2 additions for computing ${\bf n}_{jk}^E$ and 2 additions in adding the dual-face contribution ${\bf n}_{jk}^E$ to ${\bf n}_{jk}$ for an edge of each triangle. The same needs to be done for 3 edges in an element, leading to the total of 18 additions and 6 multiplications. We also need to compute the geometric centroid of an element: 4 additions and 2 multiplications (division by 3) per element. Hence, in total, it requires 22 additions and 8 multiplications per element. For a given grid, therefore, it will require $22 N_E$ additions and $8 N_E$ multiplications, where $N_E$ is the total number of triangles. 

On the other hand, the proposed algorithm described in Section \ref{eb_njk_efficient_tria} requires 2 additions and 2 multiplications for computing ${\bf n}_{j}^{E}/3$ and 2 additions for adding it to ${\bf n}_{jk} $. Assuming we compute ${\bf n}_{j}^{E}$ for each edge of a triangular element, we find it will require $12$ additions and 6 multiplications per element, resulting in $12 N_E$ additions and $6 N_E$ multiplications in total, which is less than $22 N_E$ additions and $8 N_E$ multiplications of the typical algorithm. 

Results are summarized in Table \ref{Tab.source_formulas_D}. It indicates that the proposed algorithm can be more efficient, although only slightly, by a factor $8/6 \approx 1.3$ at best in two dimensions if the cost of multiplications is dominant. 
\newline
\newline
\noindent
{\bf Remark}: 
Note that the factor $1/3$ may be applied afterwards in an edge loop after summing ${\bf n}_{j}^{E}$ as in the tetrahedral case, but it does not reduce the total cost because the total of $2 N_e$ multiplications ($N_e$ is the number of edges in a grid) is roughly equal to $6 N_E$ multiplications because we have $N_e \approx 3 N_E$ in a triangular grid. 

%Note also that the actual cost can be smaller, depending on the edge orientation: e.g., $ {\bf n}_{j}^{E_1} $ used for the edge $[j,k]$ can be used also for the edge $[j,r]$ if it is oriented from $j$ to $r$. 

%%%%%%%%%%%%%%%%%%%%%%%%%%%%%%%%%%%%%%%%%%%%%%%%%%
\begin{table}[t]
\ra{1.2}
\begin{center}
\begin{tabular}{rrrrr}\hline\hline 
\multicolumn{1}{c}{ }                                        &
\multicolumn{2}{c}{ Traditional algorithm}    &
\multicolumn{2}{c}{      Proposed algorithm}   \\  
&
Additions & Multiplications &
Additions & Multiplications 
%%% 
\\ \hline 
\\ [-0.8em]
 Triangular grid  &  $22 N_E$  &  $8 N_E$      &     $12 N_E$   & $6 N_E$  
\\ [0.8em]
%Tetrahedral grid & $153 N_E$  &  $93 N_E + 3 N_e$  &  $54 N_E$   & $36   N_E + 3 N_e$ 
Tetrahedral grid & $205 N_E$  &  $129.5 N_E $  &  $72 N_E$   & $39.5  N_E$ 
\\ [-0.8em]
  \\  \hline  \hline  
\end{tabular}
\caption{Complexity comparison for lumped directed-area vector computations in a given grid by traditional and proposed algorithms. $N_E$ is the total number of 
triangular/tetrahedral elements in a given grid. 
}
\label{Tab.source_formulas_D}
\end{center}
\end{table}
%%%%%%%%%%%%%%%%%%%%%%%%%%%%%%%%%%%%%%%%%%%%%%%%%%

%%%%%%%%%%%%%%%%%%%%%%%%%%%%%%%%%%%%%%%%%%%%%%%%
\subsubsection{Tetrahedral Grid}
\label{eb_njk_complexity_tet}
%%%%%%%%%%%%%%%%%%%%%%%%%%%%%%%%%%%%%%%%%%%%%%%%

We first consider the traditional algorithm described in Section \ref{eb_njk_typical_tet}. This particular implementation requires, at each edge of an element within a loop over elements, 
3 additions and 3 multiplications for computing $(x_m,y_m,z_m)$, 
6 additions and 3 multiplications for computing $(x_L, y_L, z_L)$, 
6 additions and 3 multiplications for computing $(x_R, y_R, z_R)$, 
6 additions and 6 multiplications for computing ${\bf n}_{jk}^{E,L}$ (the factor $1/2$ ignored here), 
6 additions and 6 multiplications for computing ${\bf n}_{jk}^{E,R}$ (the factor $1/2$ ignored here), 
3 additions for adding them, and 
3 additions for adding the result to ${\bf n}_{jk}$. 
They add up to 33 additions and 21 multiplications. All these are required for the total of 6 edges in each element, 
resulting in 196 additions and 126 multiplications. 
For each element, we also need to compute the element centroid $(x_E,y_E,z_E)$, which requires 
9 additions and 3 multiplications. 
Therefore, at each element, it requires 205 additions and 129 multiplications. Finally, in total for a tetrahedral grid, it will 
require $205 N_E$ additions and $129 N_E + 3 N_e$ multiplications, where $3 N_e$ multiplications is the cost of applying the factor $1/2$ in an edge loop, and $N_E$ and $N_e$ are the numbers of tetrahedra and edges in a grid, respectively. For a tetrahedral grid, we typically have $N_e \approx 7 N_v $ and $N_E \approx 6 N_v $, where $N_v$ is the total number of nodes in a grid, and therefore $N_e \approx 7/ 6 N_E$. Hence, we have $126 N_E + 3 N_e = 129.5 N_E$ multiplications. 

On the other hand, the proposed algorithm described in Section \ref{eb_njk_efficient_tet} will require, at each edge of an element within a loop over elements, 
9 additions and 6 multiplications for computing $2 {\bf n}_{j}^{E}$, and 
3 additions for adding $2 {\bf n}_{j}^{E}$ to ${\bf n}_{jk}$ (see Equation (\ref{add_contribution_to_njk_tetra})). 
 Assuming that we compute ${\bf n}_{j}^{E}$ for each edge of a tetrahedral element (6 edges in total), we find it will require 
 72 additions and 36 multiplications per element. Therefore, in total, it will cost $72 N_E$ additions and $36 N_E $ multiplications, with $3 N_e \approx 3.5 N_E $  multiplications for applying the factor $1/12$ in the subsequent edge loop as in Equation (\ref{apply_1over12_at_the_end_njk_tetra}). This algorithm is significantly cheaper than the traditional algorithm of $205 N_E$ additions and $129.5 N_E $ multiplications. 
 
Results are summarized in Table \ref{Tab.source_formulas_D}. It indicates approximately a factor 3 speed-up; actual numerical results show a factor between 2 and 3. 
 %As in the two-dimensional case, the actual cost can be smaller, depending on the edge orientation: e.g., $ {\bf n}_{j}^{E} $ used for the edge $[j,k]$ can be used also for the edge $[j,\ell]$ if it is oriented from $j$ to $\ell$, and even for $[j,r]$ if it is also oriented from $j$ to $r$ (see Figure \ref{fig:njk_contribution_3d}). 
\newline
\newline
\noindent
{\bf Remark}: The complexity could be reduced further. In three dimensions, the vector ${\bf n}_{jk}^E$ computed for the edge $\{ j,k \}$ can be used also for $\{ j,\ell \}$ and $\{ j, r \}$ if these two edges are oriented from $j$ to $\ell$ and $r$, respectively. If in addition, the edges $\{ k,\ell \}$ and $\{ k,r \}$ are oriented from $k$ to $\ell$ and $r$, respectively, we can use ${\bf n}_k^E$ for both. Then, we are left with the edge $\{ \ell ,r \}$, which will require ${\bf n}_\ell^E$ or ${\bf n}_r^E$. Hence, it will require three directed-area vectors per element while the estimate in Table \ref{Tab.source_formulas_D} assumes one directed-area computation per edge, thus six in total per element. It is difficult to estimate the true cost accurately because it depends on the edge orientations (it may require extra if-statements or extra memory). This implementation is not considered in this paper.

%%%%%%%%%%%%%%%%%%%%%%%%%%%%%%%%%%%%%%%%%%%%%%%%
\section{Computations of Dual Volumes}
\label{eb_vj}
%%%%%%%%%%%%%%%%%%%%%%%%%%%%%%%%%%%%%%%%%%%%%%%%
 
There are two approaches to computing dual volumes. One is based on a loop over edges and the other is based on a loop over elements. These algorithms are already known, and both do not require the construction of median dual control volumes. Here, we provide their derivations and a complexity analysis for the sake of completeness.

%%%%%%%%%%%%%%%%%%%%%%%%%%%%%%%%%%%%%%%%%%%%%%%%
\subsection{Edge-Based Algorithms}
\label{eb_vj_edge}
%%%%%%%%%%%%%%%%%%%%%%%%%%%%%%%%%%%%%%%%%%%%%%%%

 \indent 

 An edge-based formula for a dual volume computation can be derived by the divergence theorem. Consider the integration of the divergence of ${\bf x}$ over a median dual volume around a node $j$ in a triangular or tetrahedral grid, which can be evaluated as the surface integral: 
\begin{eqnarray}
 \int_{V_j} 
 \mbox{div} {\bf x} =  \oint_{\partial V_j} {\bf x} \cdot \hat{\bf n} dS, 
\end{eqnarray}
where $\hat{\bf n}$ denotes the unit outward normal vector of the infinitesimal surface area $dS$, or  
\begin{eqnarray}
 \int_{V_j} 
 \mbox{div} {\bf x} =  \oint_{\partial V_j} ( {\bf x} - {\bf x}_j ) \cdot \hat{\bf n} dS, 
 \label{div_theorem_v2}
\end{eqnarray}
since $\oint_{\partial V_j}   {\bf x}_j  \cdot \hat{\bf n} dS =  {\bf x}_j  \cdot  \oint_{\partial V_j}  \hat{\bf n} dS = 0 $. 
Then, the left hand side of Equation (\ref{div_theorem_v2}) simplifies to $D V_j$ (e.g., $ \mbox{div} {\bf x}  = \partial_x x + \partial_y y +  \partial_z z = 1+ 1  + 1 = 3$ in three dimensions) and the right hand side can be evaluated exactly by the edge-based discretization:
\begin{eqnarray}
D V_j =  \sum_{ k \in \{ k_j \}}    \Phi_{jk}  |  {\bf n}_{jk} |,
\end{eqnarray}
where the numerical flux $ \Phi_{jk} $ is given by the averaged flux,
\begin{eqnarray}
  \Phi_{jk} = 
\frac{ (  {\bf x}_j - {\bf x}_j  ) + (  {\bf x}_k - {\bf x}_j  )  }{2}  \cdot \hat{\bf n}_{jk}  =  \frac{   {\bf x}_k  -   {\bf x}_j    }{2}  \cdot \hat{\bf n}_{jk}.
\end{eqnarray}
Or we can write simply, 
\begin{eqnarray}
D V_j =  \sum_{ k \in \{ k_j \}}  \frac{   {\bf x}_k  -   {\bf x}_j    }{2}  \cdot  {\bf n}_{jk},
\end{eqnarray}
which gives
\begin{eqnarray}
  V_j =  \frac{1}{2 D} \sum_{ k \in \{ k_j \}} \left(    {\bf x}_k  -    {\bf x}_j \right)  \cdot  {\bf n}_{jk}. 
\end{eqnarray}
%Note that this formula does not require any boundary closure because $( {\bf x} - {\bf x}_j ) \cdot \hat{\bf n} = 0$ over a boundary face that is flat. 
This can be implemented in a loop over edges, distributing the dot product to the two end nodes at each edge, 
\begin{eqnarray}
  V_j =  V_j  +  \left(    {\bf x}_k  -    {\bf x}_j \right)  \cdot  {\bf n}_{jk}, 
  \quad
  V_k =  V_k +  \left(    {\bf x}_k  -    {\bf x}_j \right)  \cdot  {\bf n}_{jk},
\end{eqnarray}
and then apply the factor $1/(2D)$ in a subsequent loop over nodes:
\begin{eqnarray}
  V_j =  \frac{1}{2 D} V_j,  \,\  \mbox{for} \,\,  j= 1,2,3, \cdots, N_v. 
\end{eqnarray}

%%%%%%%%%%%%%%%%%%%%%%%%%%%%%%%%%%%%%%%%%
%  
%%%%%%%%%%%%%%%%%%%%%%%%%%%%%%%%%%%%%%%%%
 % \begin{figure}[htbp!]
  \begin{figure}[t]
    \centering
          \begin{subfigure}[t]{0.48\textwidth}
    \centering
        \includegraphics[width=0.8\textwidth,trim=0 0 0 0,clip]{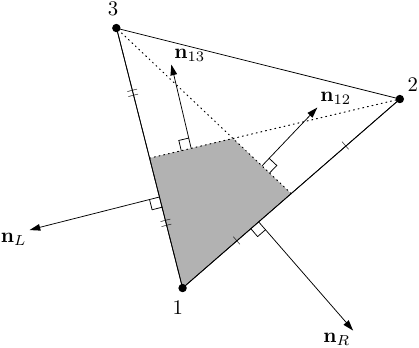}
          \caption{Triangular element.}
\label{fig:dual_contribution_2d}%
      \end{subfigure}
      \hfill
          \begin{subfigure}[t]{0.48\textwidth}
    \centering
        \includegraphics[width=0.8\textwidth,trim=0 0 0 0,clip]{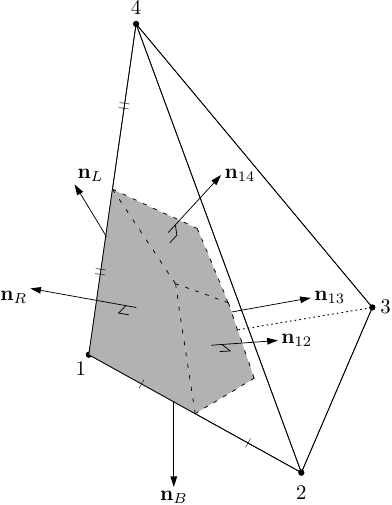}
          \caption{Tetrahedral element.}
\label{fig:dual_contribution_3d}%
      \end{subfigure}
            \caption{
\label{fig:dual_contribution}%
Median dual volume contributions to a node $1$ within an element.  } 
\end{figure}
%%%%%%%%%%%%%%%%%%%%%%%%%%%%%%%%%%%%%%%%%
% 
%%%%%%%%%%%%%%%%%%%%%%%%%%%%%%%%%%%%%%%%%

%%%%%%%%%%%%%%%%%%%%%%%%%%%%%%%%%%%%%%%%%%%%%%%%
\subsection{Element-Based Algorithms}
\label{eb_vj_elm}
%%%%%%%%%%%%%%%%%%%%%%%%%%%%%%%%%%%%%%%%%%%%%%%%

\indent 
  
Alternatively, we can compute dual volumes in a loop over elements. In fact, the dual volume contribution from a simplex element to a node is exactly $1/(D+1)$ of the element volume for $D=2$ and $D=3$. For a triangular element as illustrated in Figure \ref{fig:dual_contribution_2d}, we apply the divergence theorem to the dual volume $V_1^E$ (shaded area) contributing to the node $1$:
\begin{eqnarray}
 \int_{V_1^E} 
 \mbox{div} {\bf x} =  \oint_{\partial V_1^E} ( {\bf x} - {\bf x}_1)  \cdot \hat{\bf n} dS.
\end{eqnarray}
The left hand side is $ 2 V_1^E$ because $ \mbox{div} {\bf x} = \partial_x x + \partial_y y = 1 + 1 = 2$, and the surface integral can be evaluated exactly by using the edge-based quadrature formula derived in Appendix B of Ref.~\cite{nishikawa:AIAA2010}:
\begin{eqnarray}
2 V_1^E 
 &=& 
 \frac{   ( {\bf x}_1 - {\bf x}_1)  +  ( {\bf x}_2 - {\bf x}_1)  }{2} \cdot  {\bf n}_{12} 
 + 
 \frac{   ( {\bf x}_1 - {\bf x}_1)  +  ( {\bf x}_3 - {\bf x}_1)  }{2} \cdot  {\bf n}_{13} 
 \nonumber \\ [2ex]
 &+&
 \left[
 \frac{5}{6}  ( {\bf x}_1 - {\bf x}_1)  +  \frac{1}{6} ( {\bf x}_3 - {\bf x}_1) 
 \right]
 \cdot  {\bf n}_{L} +
 \left[
 \frac{5}{6}  ( {\bf x}_1 - {\bf x}_1)  +  \frac{1}{6}  ( {\bf x}_2 - {\bf x}_1)  
 \right]
 \cdot  {\bf n}_{R},
\end{eqnarray}
where ${\bf n}_{L}  = {\bf n}_{2} /2 $ and ${\bf n}_{R}  = {\bf n}_{3} /2 $. But the last two terms vanish by orthogonality, and we are left with
\begin{eqnarray}
2 V_1^E 
 = 
 \frac{    ( {\bf x}_2 - {\bf x}_1)  }{2} \cdot  {\bf n}_{12} 
 + 
 \frac{    ( {\bf x}_3 - {\bf x}_1)  }{2} \cdot  {\bf n}_{13}.
 \label{eq_xxx}
\end{eqnarray}
Furthermore, by the formula (\ref{njk_E1_2d}), we have
\begin{eqnarray}
 {\bf n}_{12}  =  \frac{1}{3}  \left[  {\bf n}_{1}^E + \frac{1}{2} {\bf n}_3 \right], \quad
 {\bf n}_{13}  =  \frac{1}{3}  \left[  {\bf n}_{1}^E + \frac{1}{2}  {\bf n}_2 \right], 
\end{eqnarray}
and substitute these into Equation (\ref{eq_xxx}), 
\begin{eqnarray}
2 V_1^E 
 &=& 
 \frac{    ( {\bf x}_2 - {\bf x}_1)  }{6} \cdot  {\bf n}_{1}^E 
 + 
 \frac{    ( {\bf x}_3 - {\bf x}_1)  }{6} \cdot  {\bf n}_{1}^E  \nonumber \\ [2ex]
 &=& 
 \frac{1}{3} V_E +  \frac{1}{3} V_E ,
\end{eqnarray}
where we have used the fact $( {\bf x}_2 - {\bf x}_1) \cdot  {\bf n}_3 =  ( {\bf x}_3 - {\bf x}_1) \cdot  {\bf n}_2 = 0$, and $V_E =  \frac{    ( {\bf x}_2 - {\bf x}_1)  }{2} \cdot  {\bf n}_{1}^E  =  \frac{    ( {\bf x}_3 - {\bf x}_1)  }{2} \cdot  {\bf n}_{1}^E$ is the volume (area) of the element $E$, thus leading to 
\begin{eqnarray}
 V_1^E  =  \frac{1}{3} V_E.
\end{eqnarray}

Similarly, for a tetrahedral element as illustrated in Figure \ref{fig:dual_contribution_3d}, we apply the divergence theorem to the medial dual volume $V_1^E$ (shaded volume):
\begin{eqnarray}
 \int_{V_1^E} 
 \mbox{div} {\bf x} =  \oint_{\partial V_1^E} ( {\bf x} - {\bf x}_1)  \cdot \hat{\bf n} dS.
\end{eqnarray}
The left hand side is $ 3 V_1^E$ because $ \mbox{div} {\bf x} = \partial_x x + \partial_y y + \partial_z z = 1 + 1 + 1 = 3$, and the surface integral can be evaluated exactly by using the edge-based quadrature formula derived in Appendix B of Ref.~\cite{nishikawa:AIAA2010}:
\begin{eqnarray}
3 V_1^E 
 &=& 
 \frac{   ( {\bf x}_1 - {\bf x}_1)  +  ( {\bf x}_2 - {\bf x}_1)  }{2} \cdot  {\bf n}_{12} 
 + 
 \frac{   ( {\bf x}_1 - {\bf x}_1)  +  ( {\bf x}_3 - {\bf x}_1)  }{2} \cdot  {\bf n}_{13} 
 + 
 \frac{   ( {\bf x}_1 - {\bf x}_1)  +  ( {\bf x}_4 - {\bf x}_1)  }{2} \cdot  {\bf n}_{14} 
 \nonumber \\ [2ex]
 &+&
 \left[
 \frac{6}{8}  ( {\bf x}_1 - {\bf x}_1)  +  \frac{1}{8} ( {\bf x}_3 - {\bf x}_1)   +  \frac{1}{8} ( {\bf x}_4 - {\bf x}_1) 
 \right]
 \cdot  {\bf n}_{L} 
 +
 \left[
\frac{6}{8}  ( {\bf x}_1 - {\bf x}_1)  +  \frac{1}{8} ( {\bf x}_4 - {\bf x}_1)   +  \frac{1}{8} ( {\bf x}_2 - {\bf x}_1) 
 \right]
 \cdot  {\bf n}_{R}
 \nonumber \\ [2ex]
 &+&
 \left[
\frac{6}{8}  ( {\bf x}_1 - {\bf x}_1)  +  \frac{1}{8} ( {\bf x}_2 - {\bf x}_1)   +  \frac{1}{8} ( {\bf x}_3 - {\bf x}_1) 
 \right]
 \cdot  {\bf n}_{B},
\end{eqnarray}
but the last three terms vanish by orthogonality and we are left with
\begin{eqnarray}
3 V_1^E 
 = 
 \frac{    ( {\bf x}_2 - {\bf x}_1)  }{2} \cdot  {\bf n}_{12} 
 + 
 \frac{    ( {\bf x}_3 - {\bf x}_1)  }{2} \cdot  {\bf n}_{13}
 + 
 \frac{    ( {\bf x}_4 - {\bf x}_1)  }{2} \cdot  {\bf n}_{14}.
 \label{eq_xxx_tet}
\end{eqnarray}
By the formula (\ref{njk_local_formula_3d}), we have
\begin{eqnarray}
 {\bf n}_{12}  = \frac{1}{6} 
  \left[
    {\bf n}_1^E   +   \frac{1}{2} (  {\bf n}_R  +  {\bf n}_B ) 
\right], 
\quad
 {\bf n}_{13}  = \frac{1}{6} 
  \left[
    {\bf n}_1^E   +   \frac{1}{2} (  {\bf n}_B  +  {\bf n}_L ) 
\right], 
\quad
 {\bf n}_{14}  = \frac{1}{6} 
  \left[
    {\bf n}_1^E   +   \frac{1}{2} (  {\bf n}_R  +  {\bf n}_L ) 
\right].
\end{eqnarray}
Substituting them in Equation (\ref{eq_xxx_tet}), we obtain
\begin{eqnarray}
3 V_1^E 
 = 
 \frac{    ( {\bf x}_2 - {\bf x}_1)  }{12} \cdot  {\bf n}_1^E
 + 
 \frac{    ( {\bf x}_3 - {\bf x}_1)  }{12} \cdot   {\bf n}_1^E
 + 
 \frac{    ( {\bf x}_4 - {\bf x}_1)  }{12} \cdot  {\bf n}_1^E,
 \label{eq_xxx_tet_2}
\end{eqnarray}
where all the terms involving ${\bf n}_L $, ${\bf n}_R$, and ${\bf n}_B$ have vanished by orthogonality. It gives
\begin{eqnarray}
3 V_1^E 
 = 
 \frac{1}{4} V_E +  \frac{1}{4} V_E +  \frac{1}{4} V_E, 
\end{eqnarray}
where $V_E =  \frac{    ( {\bf x}_i - {\bf x}_1)  }{3} \cdot  {\bf n}_1^E $, $i=1,2,3$, is the volume of the element $E$, and thus 
\begin{eqnarray}
 V_1^E 
 = 
 \frac{1}{4} V_E.
\end{eqnarray}

Therefore, the dual volume contribution from an element to a node is exactly $1/(D+1)$ of the element volume, where $D=2$ for a triangular grid and $D=3$ for a tetrahedral grid. Then, we can compute the dual volumes as follows: within a loop over elements, distribute the contribution to a node $i$ of an element $E$ as
\begin{eqnarray}
    V_i = V_i + \frac{1}{D +1} V_E,  \quad  i  \in \{ i_E \} ,    \,\  \mbox{for} \,\,  E= 1,2,3, \cdots, N_E, 
\end{eqnarray}
where $\{ i_E \}$ is a set of nodes of an element $E$ (3 nodes for a triangle and 4 nodes for a tetrahedron), and $V_E$ is the volume of the element $E$ (area in two dimensions). Or, in order to reduce the number of multiplications, we may apply the factor $ \frac{1}{D +1}$ in a separate loop over nodes after accumulating $V_E$: accumulate local contributions,  
\begin{eqnarray}
    V_i = V_i + V_E,  \quad  i  \in \{ i_E \}  ,   \,\  \mbox{for} \,\,  E= 1,2,3, \cdots, N_E, 
\end{eqnarray}
for all elements, and then apply the factor $ \frac{1}{D +1}$ in a loop over nodes, 
\begin{eqnarray}
  V_j = \frac{1}{D +1} V_j ,  \,\  \mbox{for} \,\,  j= 1,2,3, \cdots, N_v. 
\end{eqnarray}

%%%%%%%%%%%%%%%%%%%%%%%%%%%%%%%%%%%%%%%%%%%%%%%%
\subsection{Complexity and Discussion}
\label{eb_vj_complexity}
%%%%%%%%%%%%%%%%%%%%%%%%%%%%%%%%%%%%%%%%%%%%%%%%

As before, we provide a complexity analysis based on the number of additions and multiplications required in each algorithm and compare them, ignoring all other factors. The element-based algorithm in Section \ref{eb_vj_elm} requires 
  5 additions and 2 multiplications for computing $V_E$ (1 cross product) in two dimensions, and 
12 additions and 9 multiplications in three dimensions. Also, one multiplication per node, thus $N_v$ multiplications in total; or approximately $0.5  N_E$ multiplications in two dimensions since $N_E \approx 2 N_v$ in a triangular grid, and $0.167 N_E $ (or roughly $0.2 N_E$) multiplications in three dimensions since $N_E \approx 6 N_v$ in a tetrahedral grid. Therefore, for a given grid, it requires the total of 
$5 N_E$ additions and $2.5 N_E $ multiplications in two dimensions, and 
$12 N_E$ additions and $ 9.2 N_E$ multiplications in three dimensions. 
%$5 N_E$ additions and $2 N_E + N_v$ multiplications in two dimensions, and 
%$12 N_E$ additions and $ 9 N_E + N_v$ multiplications in three dimensions.

%%%%%%%%%%%%%%%%%%%%%%%%%%%%%%%%%%%%%%%%%%%%%%%%%%
\begin{table}[t]
\ra{1.2}
\begin{center}
\begin{tabular}{rrrrr}\hline\hline 
\multicolumn{1}{c}{ }                                        &
\multicolumn{2}{c}{ Element-based algorithm}    &
\multicolumn{2}{c}{ Edge-based algorithm}   \\  
&
Additions & Multiplications &
Additions & Multiplications 
%%% 
\\ \hline 
\\ [-0.8em]
 Triangular grid  &  $5 N_E$  &  $2.5 N_E $      &   $4.5 N_E$   & $3.5 N_E $ 
\\ [0.8em] 
Tetrahedral grid & $12 N_E$  &  $9.2 N_E $  &     $5.8 N_E$   & $3.7 N_E $ 
\\ [-0.8em]
  \\  \hline  \hline  
\end{tabular}
\caption{Complexity comparison for median dual volume computations in a given grid by element-based and edge-based algorithms. $N_E$ is the total number of 
triangular/tetrahedral elements in the grid. 
}
\label{Tab.dual_volume_D}
\end{center}
\end{table}
%%%%%%%%%%%%%%%%%%%%%%%%%%%%%%%%%%%%%%%%%%%%%%%%%%

On the other hand, the edge-based algorithm in Section \ref{eb_vj_edge} requires 
3 additions and 2 multiplications per edge in two dimensions, and 
5 additions and 3 multiplications per edge in three dimensions. Also, one multiplication per node, thus $N_v$ multiplications in total. 
Therefore, it requires the total of 
$3 N_e$ additions and $2 N_e +  N_v \approx 2 N_e + 0.5 N_E$ multiplications in two dimensions, and 
$5 N_e$ additions and $3 N_e + N_v \approx 2 N_e + 0.2 N_E$ multiplications in three dimensions. Since we have $N_e \approx 3/2 N_E$ in two dimensions and $N_e \approx 7/6 N_E$ for large grids, these estimates can be expressed approximately as 
$4.5 N_E$ additions and $3.5 N_E $ multiplications in two dimensions, and $5.8 N_E$ additions and $3.7 N_E$ multiplications in three dimensions. 

The results are summarized in Table \ref{Tab.dual_volume_D}. As can be seen, the edge-based algorithm is not necessarily more efficient than the element-based algorithm in two dimensions. However, the edge-based algorithm is expected to be more efficient than the element-based algorithm in three dimensions, at least by a factor of 2. 
 
One way to verify the dual volume data is to compute the sum of the dual volumes at all nodes and confirm that it is equal to the volume of an entire domain covered by a grid. This is typically performed by using a grid in a unit square in two dimensions and a unit cube in three dimensions: the sum of the dual volumes must be 1. 
 
%Since $ N_e \approx 7 N_v$ and $N_E \approx 6 N_v$ in a typical tetrahedral grid, where $N_v$ is the total number of nodes in a tetrahedral grid, we find that the above formula requires $21 N_v$ additions and $28 N_v$ multiplications whereas the previous approach requires $72 N_v$ additions and $ 54  N_v$ multiplication. Therefore, the above edge-based algorithm is more efficient. The edge-based algorithm has been known and utilized in the source discretization for the third-order edge-based discretization method \cite{nishikawa_liu_source_quadrature:jcp2017}. 

%   Elm based : $  8 N_E$ additions   and $ 2 N_E$ multiplications in 2D  + 1 multiplication per node.
% Edge based : $    5 N_e$ additions and $ 2 N_e$ multiplications in 2D  + 1 multiplication per node.      3*NE = 2*Nei + Neb,  ->  Nei = 1.5*NE       NE = 2*Nv -> Nv = 0.5*NE.
%                ---> $7.5 N_E$ additions  and $ 3 N_E$ multiplications in 2D.

Note that there is no need to form median dual control volumes in both algorithms. As the lumped directed-area vectors can be computed also without forming the median dual control volumes, there will be no need at all to form median dual control volumes in a code. It greatly simplifies the implementation of the edge-based discretization.

%%%%%%%%%%%%%%%%%%%%%%%%%%%%%%%%%%%%%%%%%%%%%%%%
\section{Results}
\label{eb_results}
%%%%%%%%%%%%%%%%%%%%%%%%%%%%%%%%%%%%%%%%%%%%%%%%
   
 \indent 

We implemented the proposed algorithms in two- and three-dimensional edge-based solvers, which are both serial codes, and verified their implementations by performing the verification tests mentioned earlier (results are not shown). Then, we tested the proposed algorithms by computing the lumped directed-area vectors and the dual volumes for a triangular grid in a unit square domain with 65,536 nodes and 130,050 triangles, and for a tetrahedral grid in a unit cube domain with 16,974,593 nodes and 100,663,296 tetrahedra. CPU times are measured in seconds for each computation, the computations are repeated 10 times, and CPU times are averaged over 10 executions. 
Results are summarized in Table \ref{Tab.results_njk_vj}.

For a triangular grid, the proposed algorithm is only slightly faster than the traditional algorithm by a factor 1.12, and actually the edge-based dual volume computation took longer than the element-based computation. These results are not far from the complexity analysis results, which roughly imply a factor 1.3 or less for the computation of ${\bf n}_{jk}$ (see Table \ref{Tab.source_formulas_D}) and the edge-based algorithm is expected to be slightly more expensive than the element-based algorithm for the computation of $V_j$ (see Table \ref{Tab.dual_volume_D}). Therefore, as far as the dual volume computation is concerned, the element-based algorithm may be preferred for a triangular grid. 

On the other hand, for a tetrahedral grid, the proposed algorithm provides a factor 2.16 speed-up for the computation of ${\bf n}_{jk}$ and 1.89 for the computation of $V_j$. Again, these results are not far from the complexity analysis results, which roughly indicate a factor 3 for the computation of ${\bf n}_{jk}$ (see Table \ref{Tab.source_formulas_D}) and 2 for the computation of $V_j$ (see Table \ref{Tab.dual_volume_D}). For time-accurate simulations with deforming grids, the grid metrics need to be recomputed at every physical time step, and the proposed algorithm can provide a significant speed-up. For example, over 1,000 time steps, the traditional algorithm would take $29.22 \times 1,000 = 29,220$ seconds ($\approx$ 8 hours) for the grid metric computation, and the proposed algorithm would take $13.5 \times 1,000 = 13,500$ seconds ($\approx$ 3.75 hours). 

There can be various factors contributing to the difference between the complexity analysis results and actual CPU time measurements. Especially, these computations are not so heavy and thus, for example, even an integer increment of a loop may have a significant contribution to the total CPU time. Further detailed studies on computing time is beyond the scope of this paper. The objective of this paper is to show that these computations can be performed without forming median dual control volumes and the proposed algorithm for the directed-area vectors is simpler and more efficient than a traditional algorithm. 
 
%%%%%%%%%%%%%%%%%%%%%%%%%%%%%%%%%%%%%%%%%%%%%%%%%%
\begin{table}[h!]
\ra{1.2}
\begin{center}
\begin{tabular}{rrrrr}\hline\hline 
\multicolumn{1}{c}{ }                                        &
\multicolumn{2}{c}{ Conventional algorithm}    &
\multicolumn{2}{c}{      Proposed algorithm}   \\  
&
${\bf n}_{jk}$ & $V_j$ &
${\bf n}_{jk}$ & $V_j$ 
%%% 
\\   \hline 
\\ [-0.5em]
 Triangular grid  &  4.09e-03 (1.12)  &  7.25e-04 (0.81)      &     3.85e-03   & 9.56e-04 
\\ [0.8em] 
Tetrahedral grid & 29.22 (2.16)  &  2.38 (1.89)  &  13.5   & 1.26
\\ [-0.8em]
  \\  \hline  \hline  
\end{tabular}
\caption{ CPU time (seconds) required to compute the lumped directed-area vectors, ${\bf n}_{jk}$, and median dual volumes, $V_j$, for a triangular grid with 65,536 nodes and 130,050 triangles and a tetrahedral grid with 16,974,593 nodes and 100,663,296 tetrahedra. Numbers in parentheses indicate the ratio to the corresponding time in the traditional algorithm (e.g.,  2.16 = 29.22/13.5). 
CPU times are averaged over 10 executions. 
}
\label{Tab.results_njk_vj}
\end{center}
\end{table}
%%%%%%%%%%%%%%%%%%%%%%%%%%%%%%%%%%%%%%%%%%%%%%%%%%

%%%%%%%%%%%%%%%%%%%%%%%%%%%%%%%%%%%%%%%%%%%%%%%%
\section{Concluding Remarks}
\label{conclusions}
%%%%%%%%%%%%%%%%%%%%%%%%%%%%%%%%%%%%%%%%%%%%%%%%

 \indent 

We have shown that the lumped directed-area vectors and the dual volumes required to implement the edge-based discretization can be computed without explicitly forming median dual control volumes and the resulting algorithms for computing the lumped directed-area vectors can be more efficient than a traditional algorithm. Numerical results show that the speed-up factor can be 1.12 for a triangular grid and 2.16 for a tetrahedral grid. More importantly, the proposed algorithm is expected to dramatically simplify the implementation of the edge-based discretization in four and higher dimensions as it does not require the construction of median dual volumes. 

Future work includes applications to practical CFD codes, especially for time-dependent simulations with deforming grids. Extensions to mixed-element grids might be possible but the edge-based discretization loses second-order accuracy unless it satisfies certain regularity conditions \cite{nishikawa:AIAA2010,Boris_Jim_NIA2007-08} or a flux correction technique \cite{Nishikawa_FC_Polyhedral:jcp2022} is incorporated. Extensions to another dual control volume defined with the face-area-weighted centroid \cite{Nishikawa_aiaa2020-1786} instead of the geometric centroid, if possible, would be useful because it is known to reduce effects of grid skewness and improve accuracy and efficiency of an edge-based solver as demonstrated in Ref.~\cite{Nishikawa_aiaa2020-1786}. Finally, extensions to four dimensions are currently underway, and will be reported in a separate paper. 

% It remains to demonstrate that the proposed algorithm is directly applicable to simplex-element grids in four dimensions. 

%%%%%%%%%%%%%%%%%%%%%%%%%%%%%%%%%%%%%%%%
\section*{Acknowledgments}
%%%%%%%%%%%%%%%%%%%%%%%%%%%%%%%%%%%%%%%%

The author gratefully acknowledges support from the U.S. Army Research Office under the contract/grant number W911NF-19-1-0429. 
%and support by Software CRADLE, part of Hexagon.

%%%%%%%%%%%%%%%%%%%%%%%%%%%%%%%%%%%%%%%%%%%%%%%%%%%%%%%%%%%
\bibliography{../../../bibtex_nishikawa_database}
\bibliographystyle{unsrt}
%%%%%%%%%%%%%%%%%%%%%%%%%%%%%%%%%%%%%%%%%%%%%%%%%%%%%%%%%%%

\end{document}